\documentclass[journal]{IEEEtran}
\usepackage{amsmath,amsfonts}
\usepackage{algorithm}
\usepackage{algpseudocode}
\usepackage{array}
\usepackage[caption=false,font=normalsize,labelfont=sf,textfont=sf]{subfig}
\usepackage{textcomp}
\usepackage{stfloats}
\usepackage{url}
\usepackage{verbatim}
\usepackage{graphicx}
\usepackage{subfig}
\usepackage{siunitx} 
\def\n2#1{\ensuremath{\SI[scientific-notation=true, round-mode=places, round-precision=2]{#1}{}}}
\usepackage[table]{xcolor}
\usepackage{xcolor, soul}
\usepackage{amsmath}
\DeclareMathOperator*{\argmax}{arg\,max}

\def\P{\mathbb{P}}
\begin{document}

\title{BOP-Elites, a Bayesian Optimisation Approach to\\ Quality Diversity Search\\ with Black-Box descriptor functions}

\author{P. Kent,
 Adam Gaier, Jean-Baptiste Mouret and Juergen Branke}

\maketitle

\begin{abstract}
Quality Diversity (QD) algorithms such as MAP-Elites are a class of optimisation techniques that attempt to find many high performing points that all behave differently according to a user-defined behavioural metric.

In this paper we propose the Bayesian Optimisation of Elites (BOP-Elites) algorithm. Designed for problems with expensive black-box fitness and behaviour functions, it is able to return a QD solution-set with excellent final performance already after a relatively small number of samples. BOP-Elites models both fitness and behavioural descriptors with Gaussian Process (GP) surrogate models and uses Bayesian Optimisation (BO) strategies for choosing points to evaluate in order to solve the quality-diversity problem. 
In addition, BOP-Elites produces high quality surrogate models which can be used after convergence to predict solutions with any behaviour in a continuous range.

An empirical comparison shows that BOP-Elites significantly outperforms other state-of-the-art algorithms without the need for problem-specific parameter tuning.
\end{abstract}


\section{Introduction}
Optimisation algorithms are ubiquitous in science, engineering, and research and typically attempt to find the best feasible solution for a given objective function. In recent years, a new class of optimisation problems has emerged, called Quality Diversity optimisation (QDO).
 
In QDO settings, alongside a solution's objective performance, solutions exhibit some behavioural characteristics or \emph{descriptor values} which maps solutions on to a low dimensional \emph{descriptor space}. In general, the idea is that these descriptors are of interest to a decision maker, providing an informative way to distinguish between solutions. By requiring a search algorithm to return points well distributed in descriptor space, we develop a solution set that is meaningfully diverse and provides insights about the best possible objective performance depending on the point in descriptor space. This process of gaining information on the relationship between an objective and its descriptors is referred to as \emph{illumination} in QD literature \cite{mouret2015illuminating}. 

One compelling argument for incorporating QD search in design tasks is that it enables the decision-maker to leverage their expertise when selecting the final solution. This approach makes it possible to consider criteria that are difficult to encode in an objective function, such as aesthetics or manufacturing ease. By using a QD search, the decision-maker can find a set of high-quality starting designs that exhibit varied characteristics, allowing them to refine the design manually while taking into account their expert knowledge of additional factors that are challenging to quantify.

Interest in QD research has rapidly grown since the introduction of the Multi-Dimensional Archive of Elites (MAP-Elites) algorithm \cite{mouret2015illuminating} that has wide reaching applications in high dimensional domains and exhibits excellent performance with relatively low computational cost. Whilst MAP-Elites is computationally efficient, it is not designed for sample efficiency, making it unsuitable for use when the objective functions are expensive. 

The Surrogate Assisted Illumination (SAIL) algorithm \cite{gaier2018data}  drastically improves sample efficiency by using a surrogate model of the objective function. In contrast to MAP-Elites, which provides its final archive of observed points as a solution set, SAIL offers a prediction based on its surrogate model. While SAIL learns a surrogate model solely for the objective, Surrogate-Assisted Phenotypic Niching (SPHEN) \cite{sphen} expands upon SAIL by modeling expensive features using a GP, thus making it suitable for use with black-box features.

In this study, we present BOP-Elites, a Bayesian optimization algorithm tailored for Quality-Diversity problems\footnote{A preliminary paper outlining the basic idea for BOP-Elites was released on arxiv in 2020 by the authors in \cite{BOPelites}}. Like SPHEN, BOP-Elites employs Gaussian process surrogate models to represent both objective and descriptor landscapes. However, in contrast to prior research, our method applies a traditional Bayesian Optimization framework by continuously optimizing an acquisition function using sequential point selection, which enables a principled and efficient search strategy with unparalleled sample efficiency. Furthermore, BOP-Elites utilises a structured archive for retaining high-performing observations and ultimately delivers a solution set of observed points, akin to the MAP-Elites method.

Inspired by the SAIL methodology, we additionally explore the 'illumination' capabilities of BOP-Elites by evaluating the quality of the surrogate models at the end of the run. We further show that surrogate models built by BOP-Elites can be used to 'upscale' a solution archive, i.e., use the surrogates trained when building a lower resolution archive, with fewer regions, to fill a higher resolution archive with high performing predictions.

The key contributions of this paper are as follows:

\begin{itemize}
    \item The first principled Bayesian Optimisation algorithm for Quality-Diversity Optimisation with expensive black box objective and descriptor functions\footnote{'While surrogate-assisted methods exist that implement BO acquisition functions, the point selection mechanism used by these algorithms deviates from BO methodology'}.
    \item An extension of the well known Expected Improvement \cite{Jones98} acquisition function for the QD case and, more generally, to the problem of optimising sub-regions of a search-space modelled with a single Gaussian Process.
    \item A method for 'upscaling' is proposed, using the surrogate models built while solving a QD problem with few regions for approximating a solution to a QD problem with many regions.
    \item A thorough experimental comparison is provided, showing BOP-Elites performance against the state-of-the-art QDO algorithms.
\end{itemize}

\section{Problem definition\label{sec:probdef}}


Given a box-constrained, $d$ dimensional search domain $X \subset {\rm I\!R}^d$ and an objective function over the search domain which returns real-valued objective and box constrained descriptor values.

\[f(x)=y, b: X   \to  {\rm I\!R\times\rm I\!R^m}.\]



The value of $y$ is the performance of the solution while $b$ contains $m$ scalar values which describe some high level characteristic in an $m$ dimensional descriptor space.  

While descriptor functions themselves need not be bounded in principle, bounds mark areas of behavioural interest. For each descriptor $b_i$, human decision makers specify a-priori a set of $N_{i}$ uniform partitions within the dimension of each descriptor $i$.  By subdividing each descriptor dimension, we form $|R| = \prod_{i=1}^m N_{i}$ regions in descriptor space with distinct combinations of behavioural qualities. A simple example would be having 2 descriptor functions each with two partitions. Suppose we wish to manufacture a part made with a combination plastic polymer, our search domain is the space of different plastics, the performance metric is the life-span of the part under use and the descriptors are observed both during manufacturing and after testing. $b_1$ is the color of the resultant polymer and has two region labels 'Mostly Blue' and 'Mostly Red'. $b_2$ measures the electrical conductivity of the part and has two labels 'Low conductivity'  and 'High conductivity'. This leads to $|R|=4$ regions with the interpretable labels $R_{1}$: 'Mostly Blue, High Conductivity', $R_{2}$: 'Mostly Blue, Low Conductivity', $R_{3}$: 'Mostly Red, High Conductivity' or $R_{4}$: 'Mostly Red, Low conductivity'.

\begin{figure}
\centering
\includegraphics[width = 0.3\textwidth]{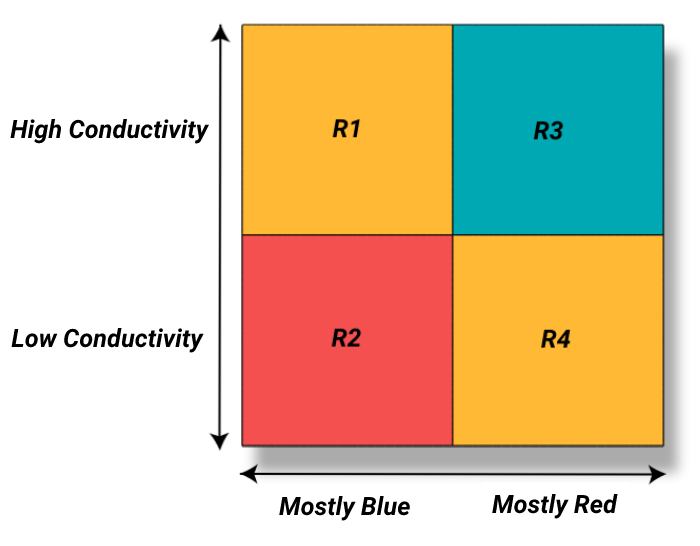}
\caption{example 2x2 solution archive for a polymer manufacturing process with 4 regions in descriptor space. Colours indicate difference in performance.}
\end{figure}

To highlight the distinctions between Quality-Diversity and other optimisation methods, it is essential to recognize that there is no universally 'better' color, as preferences vary depending on the decision-maker. Similarly, a manufacturer may desire a material with either conductive or insulating properties, depending on its intended application. The descriptors here are not treated as objectives to be optimized, but dimensions for generating diverse solutions over.

In this paper, we explore two types of problems distinguished by the characteristic of the descriptor functions. In the simpler case, the descriptor functions are  analytically known or cheap to evaluate. In the more challenging case, the descriptor functions are unknown and potentially expensive to evaluate. These descriptor functions are sometimes observable independent of the objective of a solution, e.g.\ we can easily calculate the area of a shape before testing its aerodynamic performance. Or they may only be observed at the same time as the objective evaluation, e.g.\ the amount of drag a shape exhibits when tested.

The goal of BOP-Elites is to find the set $S^*$ of optimal solutions within each achievable predefined region in descriptor space in as sample efficient a way as possible. Formally this can be expressed as
\[S^* = \{e_1^*,...,e_{|R|}^*\}\]
\[e_r^*= \argmax_x f(x) \,\,\,\, s.t \,\,\,\, x \in X,\,\, b \rightarrow R_r,  \]  

where by ``$b\rightarrow R_r$" we mean that the descriptor values $b$ place the solution in region $R_r$.
When the budget is sufficient to evaluate a point in each region, a final solution will consist of a set of `elites' $S = \{e_{1},...,e_{|R|}\}$, the best solution found for each region. In this paper, solution sets are evaluated using the Quality Diversity score (QDs) \cite{mouret2015illuminating}, which is the sum of performances of the obtained observations, 
\begin{equation}
     QDs=\sum_{r=1}^{|R|} f(e_r). \label{TE}
\end{equation} 

As we assume our domains are deterministic this can be considered a risk free solution set as all recommended solutions have been observed and their values and descriptors are known. 

If the number of regions is large and the evaluation budget small, $e_i$ may be empty for some $i$. In this case, the algorithm would return some minimum value, often 0, for all unobserved regions and the overall performance would be poor. An alternative approach, used in surrogate-assisted QD methods is to build a representative model of the objective function and attempt to predict the performance with a so-called predicted map \cite{gaier2018data}. As BOP-Elites is a model-based Bayesian Optimisation algorithm, it builds a surrogate model which can be used to predict the elite points $\hat{e}$ and the quality of such a predicted map is the true QDs value of the predicted points:

\begin{equation}
     QDs_p=\sum_{r=1}^{|R|} f(\hat{e_r}) \label{TEP}
\end{equation} 
\begin{equation*}
f( \hat{e}_r ) = \begin{cases}
f( \hat{e}_r ) &\text{if $b(x) \rightarrow R_r$}\\
0 &\text{otherwise}
\end{cases} 
\end{equation*}

In BOP-Elites, black-box descriptor functions can be predicted using the descriptor GP models. When descriptor values are mispredicted, i.e., a solution is proposed with one predicted behavioural quality which, when evaluated, is different from it's expectation, we simply set the value of a solution to 0, or a contextually meaningful lower bound.

\section{Related work\label{sec:related}}
\subsection {Quality-Diversity Algorithms}

QD has its origins in the field of evolutionary computation where, when searching for a single global optimum, maintaining genetic diversity is important to avoid premature convergence and to escape local optima \cite{ursem2002diversity}. In QDO, however, we are not seeking diversity over the search space but searching for diversity in descriptor space. In other words, QDO wishes to solve many different optimisation problems simultaneously, where each problem is defined over a shared problem space and solves the problem in a behaviourally diverse way.

By contrast to single objective optimisation, QD algorithms such as Novelty Search with local competition (NSLC) \cite{lehman2011evolving} and MAP-Elites \cite{mouret2015illuminating} return a \textit{set} of high performing solutions. The solution set returned by MAP-Elites is an example of a structured archive, a grid of points where each bin in the grid stores a solution for a region in descriptor space. By contrast NSLC returns an unstructured archive, a population of points of fixed size that sequentially compete, being replaced by points that improve diversity or performance of the population.

The QD problem is similar, but different to other classes of optimisation problems where evolutionary algorithms have been successful:
\begin{itemize}
    \item In multimodal optimisation, one is looking for several local optima in the search space, but there is only the objective function (no descriptors that define separate partitions, the partitions are implicitly defined by the location of local optima).
    \item In multi-objective optimisation, one is looking for several solutions with different trade-offs in objectives. Objectives have to be optimised simultaneously, but there are no pre-defined partitions in descriptor space. 
    \item In multi-task optimisation, several related optimisation problems have to be solved simultaneously. Our problem setting could be framed as multi-task optimisation, where each task is to find the best solution belonging to a particular region of descriptor space, i.e., all tasks have the same objective and search space but different complex constraints (descriptor regions).
    \item In constrained optimisation, there are usually just two categories, feasible and infeasible solutions, and we would like to find the best feasible solution. But identifying the best solution in a particular region of descriptor space could be seen as a constrained optimisation problem, where the region boundaries are the constraints.
\end{itemize}

QD is a vibrant field of research with recent applications in game level design, \cite{fontaine2020illuminating}, urban design, \cite{galanos2021arch}, constrained optimisation \cite{fioravanzomap}, robotics \cite{cully2015robots,nilsson2021policy,paolo2021sparse,nordmoen2021map} and the design of nanomaterials \cite{jiang2022artificial}. There are many real world situations where producing a diverse set of high performing solutions may be desirable and with such rapid development in the field \cite{fontaine2021differentiable}, the range of applications for such algorithms is likely to increase. 

A variety of algorithms have been proposed for QD search, and \cite{QDreview} provides a good overview. In the following, the focus is on the algorithms most relevant to our work, namely MAP-Elites, SAIL, and SPHEN.

{

\begin{table}
\centering
\caption{Comparison of features of QD algorithms compared in this paper, High dimensionality implies $<20$ parameters.}
\begin{tabular}{|l|c|c|c|c|}
\hline
\multicolumn{1}{|c|}{\textbf{Use-Cases}} & Map-Elites & SAIL & SPHEN & BOP-Elites \\ \hline
Blackbox-Descriptors  &      \checkmark      &     & \checkmark & \checkmark            \\ 
High dimensional $X$    &  \checkmark          &      &&            \\ \hline
\multicolumn{1}{|c|}{\textbf{Output}} & Map-Elites & SAIL & SPHEN & BOP-Elites \\ \hline
Observed archive &        \checkmark     &   \checkmark   & \checkmark & \checkmark            \\
Prediction Map &            & \checkmark     & \checkmark  & \checkmark            \\ \hline
\end{tabular} \label{tab:algcomp}
\end{table}

}
\subsection{MAP-Elites\label{sec:map-elites}}
Most modern QD takes place in the context of the MAP-Elites framework \cite{mouret2015illuminating}. An initial set of points are evaluated on the objective and descriptor functions, the best point in each region is considered an 'elite' and is stored in the archive. Elite points are randomly selected as parents from the archive and combined to produce children that are then evaluated. If a child outperforms the elite for its region in the archive, or is the first individual in that region, it becomes the new elite in that region. This relatively simple, yet powerful approach provides good coverage over the regions and finds high-performing points when given a sufficient budget.

\begin{figure}
    \centering
    \includegraphics[width = 0.40\textwidth]{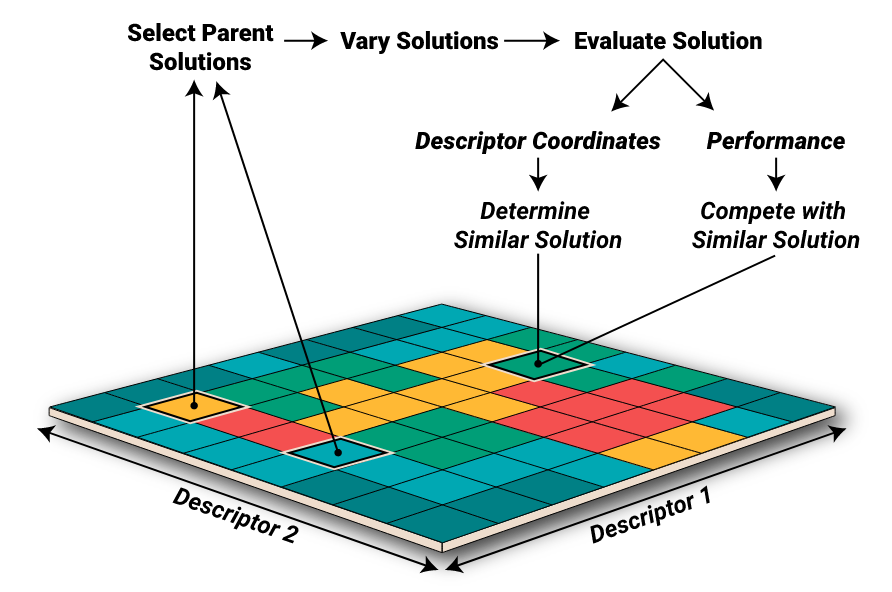}
    \caption{A visualisation of Map-Elites, image adapted from \cite{mouret2020evolving}}
    \label{fig:mapelites}
\end{figure}

\subsection{Bayesian Optimisation\label{sec:bayesopt}}

Bayesian Optimisation (BO) is an optimisation approach for expensive black-box functions, see, e.g., \cite{frazier2018tutorial} for an introduction. The general idea of BO is to start with some initial design points obtained e.g.\ by Latin Hypercube Sampling, then build a surrogate model of the objective function that not only allows to predict the quality of any solution candidate, but also provides a confidence measure for this prediction. This information is then iteratively used to decide which solution candidate to evaluate next, and update the surrogate model. Point selection is performed by optimising a so-called  \emph{acquisition function}, a mathematically motivated heuristic that balances exploration and exploitation. 
While this means that an inner optimisation problem has to be solved in every iteration, this is only based on the quick-to-evaluate acquisition function, whereas BO exhibits a very good sample efficiency in terms of the often very expensive objective function.

The most widely used BO surrogate model, and the one implemented in BOP-Elites, is the Gaussian Process (GP) model \cite{RassmussenWilliams}. Assuming that the latent function $f(.)$ can be suitably modelled with a GP, for any finite set of observed data pairs $\tau^n= \{(x_1,y_1),...,(x_n,y_n)\}$, $X_n = \{x_1,...,x_n\}, Y_n = \{y_1,...,y_n\}$  this function is modelled as a multivariate Gaussian random variable. Assuming a prior mean function $m(x) $,  the function value $f(x)$ at an unobserved point can be predicted with posterior mean and variance 

\begin{equation}     
    {\rm I\!E} \left[ f(x)| \tau^n \right]= \mu(x) =m(x)  + K_*K^{-1}(Y-m(x)) \label{GP2}
\end{equation} 
\begin{equation}
    var(f(x)) = s_0^2 = K^\top_{**}-K_*K^{-1}K_* \label{GP3}
\end{equation} 

where $K$ is the kernel, or covariance, function. $K_* = K(x,X_n)$, $K^{-1} = K(X_n,X_n)^{-1}$ and $K_{**}=K(x,x)$. While many kernels exist, throughout this paper we use the Matérn 5/2 kernel kernel in the BOP-Elites algorithm:

\begin{equation}
k(x, x') = \sigma_f^2 \sum_{d=1}^D\left[\left(1 + \frac{\sqrt{5} r}{\ell_d} + \frac{5 r^2}{3 \ell_d^2}\right) \exp\left(-\frac{\sqrt{5} r}{\ell_d}\right)\right]
\end{equation}

The kernel has $D+2$ free parameters which must be estimated from the data; $D$ characteristic length-scales $\ell_d $ which encodes the importance of the correlation in dimension $d$ for the similarity in $f$ and signal variance $\sigma_f^2 $ which sets the maximum covariance for the process. 

This derivation is for the modelling of deterministic smooth functions and therefore does not include a meaningful noise term but an adaptation to noisy functions would be straightforward. GP regression is covered in depth in the influential book by \cite{RassmussenWilliams}.

Key to BO algorithms is the acquisition function (sometimes called infill criterion), a calculation performed on the posterior distribution, which attempts to predict the value of sampling new points in the input space. Acquisition functions are problem specific and designed to balance exploration and exploitation in the search. Examples include upper confidence bound (UCB), probability of improvement (PI) \cite{Kushner64}, Expected Improvement (EI) \cite{mockus1978application,Jones98} and Knowledge Gradient (KG)\cite{Scott11}. In the following we present UCB and EI in more detail, as they will be used later on.


The UCB algorithm was first developed for the multi-armed bandit problem and implemented as a Bayesian acquisition function in \cite{srinivas2009gaussian}. It attempts to minimise regret by choosing to value points as the upper limit of a confidence bound using the posterior mean $\mu(x)$ and standard deviation $s_0$ defined as 

\begin{eqnarray}
 UCB(x) &=&\mu(x)  + \beta s_0. \label{UCB}
\end{eqnarray}

UCB algorithms have a parameter $\beta$ which is typically tuned to the problem of interest. $\beta$ controls how highly to value reducing posterior uncertainty and therefore can be considered an exploration parameter.  


Expected Improvement (EI), which we adapt for BOP-Elites, is a very popular acquisition function as it is both powerful and relatively simple to implement. 

Consider some candidate point $x$ and $\mu(x)$ the posterior predicted valuation of $x$ from our GP model. EI chooses the point for sampling that will maximise the improvement over our current best performing observation $\{x^*,y^*\}$, i.e., the point $x$ that maximizes

\begin{eqnarray}
 EI(x) &=& \mathbb{E}\left[\max(f(x)-y^*,0)\right]. \label{EI1}
\end{eqnarray}
This may be calculated in closed form \cite{Jones98,mockus1978application}
\begin{align}
 EI(x) = (\mu_0-y^*)\Phi\left(\frac{\mu_0-y^*}{s_0}\right)+s_0\phi\left(\frac{\mu_0-y^*}{s_0}\right),
\end{align}

where $\Phi$ and $\phi$ are the standard cumulative normal and density functions, $\mu_0$ and $s_0$ are the posterior mean prediction and standard deviation at $x$ from the GP. EI suggests points that have a high probability of yielding improvement either because the surrogate model predicts an improvement with high certainty or because a point has high uncertainty that could yield a high positive value. 

\subsection{Surrogate Assisted Illumination (SAIL) algorithm\label{sec:SAIL}}

Surrogate-assisted approaches have shown promise in accelerating problems which are both high-dimensional and expensive~\cite{dsame_gecco,dsame_iclr}. 

The Surrogate Assisted Illumination algorithm (SAIL) \cite{gaier2018data} is a QD algorithm that attempts to leverage the power of a GP surrogate model of the objective landscape. SAIL performs MAP-Elites using Upper Confidence Bound (UCB) as an objective function on the surrogate model. The underlying idea is to attempt to build a good model of the objective function, then  use it to predict elites. Once complete, SAIL outputs a \emph{predicted} elite for each region, and a model of the objective function which may provide additional insights. SAIL does not model the descriptor space and instead assumes that solution descriptors can be quickly and accurately obtained.

SAIL was the first approach to use surrogate models and ideas from  Bayesian optimisation  to enhance the performance and increase sample efficiency of MAP-Elites  \cite{sail_gecco}. While MAP-Elites is effective, it typically takes many observations to reach a high performing solution.  SAIL uses the evolutionary algorithm approach of MAP-Elites on the cheap-to-evaluate surrogate model, with the UCB function evlauated on the GP as the objective function. This approach explores points with high objective performance and high variance and generates both a set of predicted elites and a GP that \textit{illuminates} the relationship between performance and behaviour in a meaningful way. SAIL has been primarily applied as a tool for exploration in search spaces where simulations or experiments are expensive, such as aerodynamic design and materials science~\cite{sail_ecj,sail_aiaa,sail_gecco,produqd,sphen,sen2021integration}. Table \ref{tab:algcomp} compares the use-cases and outputs from the various algorithms considered in this paper. 

\subsection{Surrogate-assisted Phenotypic Niching}\label{SPHEN}
The general framework introduced by SAIL for using surrogate models to assist QDO has been adapted to predict descriptors along with objectives~\cite{sphen} in order to find solutions with specific features and as a component in a larger generative design process that is iteratively applied to narrow in on promising design spaces~\cite{produqd,zhang2022deep}. 

SPHEN builds upon SAIL by modelling the descriptor space with a seperate GP surrogate model. SPHEN uses the posterior mean of the GP to predict the region in descriptor space a solution is likely belonging to, and accepts the prediction taking no account of the uncertainty in the model. This allows the algorithm to predict region membership and implement SAIL where the descriptors are black box. SPHEN does not directly drive exploration in the descriptor/phenotype space, instead it updates its descriptor model indirectly through the observations gained using SAILs acquisition function over the objective.

\subsection{Bayesian Optimisation for Quality-Diversity}
While the SAIL algorithm utilises UCB, a BO acquisition function, to search for the best points to sample, its selection criterion departs from traditional BO methodology. In the acquisition step SAIL fills an \textit{acquisition map} with points that score highly on the UCB function before randomly selecting points from this acquisition map to evaluate on the true objective function. This selection process induces the possibility of randomly selecting a point from a region for which there is no real improvement to be found. Traditional EI, by contrast, attempts to compare the improvement a point provides to the objective function and selects the point for which there is the highest improvement.

In practice EI often outperforms UCB in direct comparisons \cite{Practicalbayesopt,greedisgood} and while the performance of the UCB algorithm can be affected by parameter selection, EI has no parameters to tune. These arguments taken together provide compelling justification for creating a principled implementation of EI for QD.

\section{BOP-Elites Algorithm\label{sec:bopelites}}

BOP-Elites is a model-based Bayesian optimisation algorithm for QD problems. It is designed to work with an expensive black-box objective function $f(x)$ and descriptor functions $b(x)$ which we assume can be suitably modelled with a GP.

\begin{figure}
    \centering
    \includegraphics[width = 1\linewidth]{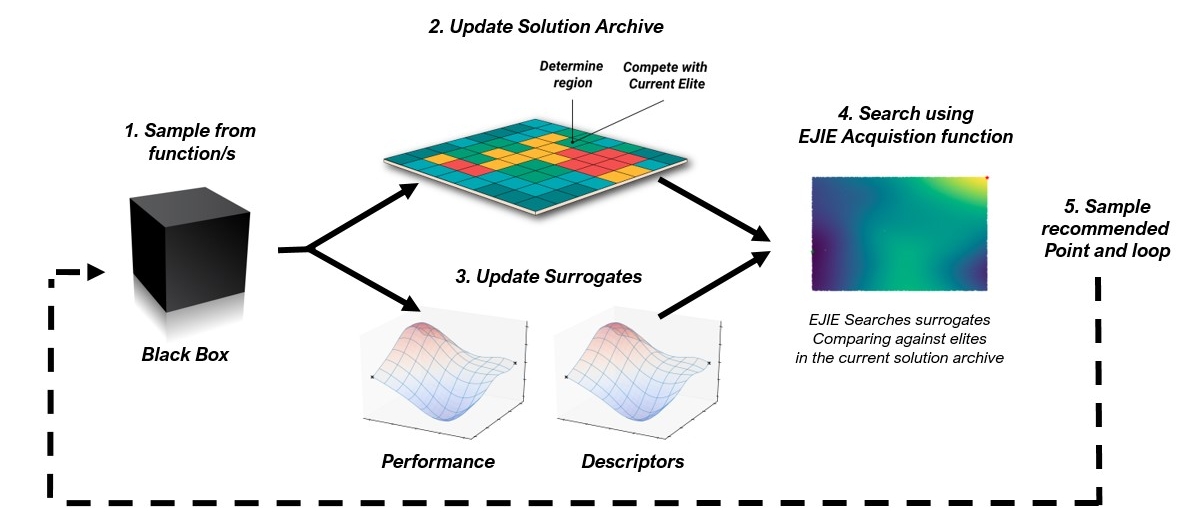}
    \caption{A visualisation of the BOP-Elites algorithm}
    \label{fig:BopAlgo}
\end{figure}

Following the Bayesian Optimisation paradigm, BOP-Elites builds GP surrogate models of $y$ and, when the descriptor function is also an expensive black-box function, each $b_i$. At each iteration, it will use its acquisition function to search the model for the next point to be sampled. After each  black-box function evaluation,  we update our models and attempt to increase the value of our solution set. The algorithm proceeds until the search budget has been exhausted. The general flow of the algorithm is visualised in Fig. \ref{fig:BopAlgo}.

BOP-Elites is initialised with a set of $n_0$ points Sobol sampled from $X$. The objective and descriptor functions are evaluated at each initial point and stored as a list of observations, the data $\mathcal{D}=\{x_i,y_i,b_i\}_{i=1}^n$. We store `elites' $\hat{S} = \{e_1,..,e_{|R|}\}$ in an archive, as in MAP-Elites, this is updated whenever a new point is sampled that performs better than the existing elite in the corresponding region.

\subsection{BOP-Elites Acquisition Function\label{acqfun}}
\subsubsection{Predicting region membership}
We propose an extension to the expected improvement acquisition function tailored to quality-diversity search, because EI is known to be effective, sample efficient and has no parameters to tune. The simplest approach would be to calculate the expected improvement for each region independently and take the highest one. However, determining which region a specific point belongs to may be computationally expensive yet is required to evaluate EI.


This issue is critical since region membership --- and hence, which elite we should compare against ---  must be predicted and comparison to the wrong elite will mis-value a point. BOP-Elites attempts to solve this issue by building surrogate models of the descriptor space and calculating the posterior probability $\mathbb{P}(x \in R_r|\mathcal{D})$ of a point belonging to region $R_r$.  

As behavioural descriptors are evaluated by $m$ descriptor functions we fit $m$ GPs using our observed data and make $m$ predictions using the posterior mean positions $\hat{b}(x) = \{ \hat{b}_1(x), ...,  \hat{b}_m(x) \}$ and variances $s(x) = \{ {s}_1(x), ...,  {s}_m(x) \}$ along each of the descriptor dimensions using Eqns.~\ref{GP2} and \ref{GP3}. Regions are formed by intersecting partitions across descriptor dimensions and the probability of a given point $x$ lying within a partition $\mathcal{P}_{j,r}$ in dimension $j$ that forms part of region $r$ is calculated as:

\begin{equation} \label{nicheprob}
\P(b_j(x) \rightarrow \mathcal{P}_{j,r}) = \Phi\left(\frac{\hat{b}_j(x)-UB}{ s_j(x)}\right)-\Phi\left(\frac{\hat{b}_j(x)-LB}{ s_j(x)}\right),
\end{equation}



where $\Phi$ is the posterior CDF of a uni-variate Gaussian and $LB$ and $UB$ are the lower and upper bounds for partition $\mathcal{P}_{j,r}$. Therefore, the likelihood of belonging to any given region is determined as the product of the individual partition probabilities. Assuming independence of the descriptors, we get the following:

\begin{equation}
    \P(x  \rightarrow R_r|\mathcal{D}) = \prod_{j}^m \P(b_j(x)  \rightarrow \mathcal{P}_{j,r}) \label{probr}
\end{equation}

\subsubsection{Calculating Expected Improvement when there is no elite}

A second issue implementing EI \cite{Jones98} arises as calculating the EI value requires comparing our posterior mean performance of a point $\mu_0$ to the value of the current elite $f(\hat{e}_r)$ in our solution archive. 

\begin{eqnarray}
    EI_r(x) &=& \mathbb{E}\left[\max(\mu_0-f(\hat{e}_r),0)\right]  \label{EIR}\\
 &=& (\mu_0-f(\hat{e_r}))\Phi\left(\frac{\mu_0-f(\hat{e_r})}{s_0}\right) \nonumber \\ & &+ \,s_0\phi\left(\frac{\mu_0-f(\hat{e_r})}{s_0}\right)   
 \label{nicheimp}
\end{eqnarray}

We propose to compare all regions, including those without an observed point, using our acquisition function, in order to identify the point that offers the best single-step improvement. Given that an empty region contributes nothing to the QD score of our archive, it seems sensible to assume that the elite performance of an empty region is zero. This means that BOP-Elites will value empty regions highly and this will drive region-filling behaviour as BOP-Elites is an improvement-greedy strategy.

\subsubsection{Expected Joint Improvement of Elites}

Combining the above, we can calculate the region-specific Expected Improvement which returns the expected improvement over the elite in region $r$ weighted by the posterior probability that it is in competition with that elite.

The proposed Expected Joint Improvement of Elites (EJIE) acquisition function takes the sum of these values for all regions to predict the improvement the point $x$ will provide to our current solution archive.

\begin{equation}
EJIE(x) = \sum_{i=1}^{R}\mathbb{P}(x \rightarrow R_r|\mathcal{D})EI_{r_i}(x) \label{EJIE}.
\end{equation}

This is  similar to the approach proposed for constrained Bayesian Optimisation (CBO) in \cite{gardner2014bayesian} where the region boundaries are constraints. The key difference is that CBO is looking for a single optimal value given some constraints while QD attempts to identify the optimal value for each of a set of non-overlapping constrained regions.

The EJIE acquisition function, Eqn.~\ref{EJIE}, gives high valuations to points that are likely to provide large improvements to the QD score. The predefined regions force the algorithm to search in areas of the input domain that single-objective optimisation would ignore, naturally diversifying the solution set. The posterior models provided by BOP-Elites furthermore provide insightful illumination of the search space, offering the opportunity for prediction of behavioural quality and objective performance at any point in the input domain with quantifiable confidence bounds. 

In the following sub-section we address difficulties with the EJIE acquisition function for QD problems and propose further adaptations to the algorithm.

\subsection{Optimisation of the acquisition function} \label{optimacq}

Optimisation of the acquisition function is performed using  pattern search (or direct search) \cite{patternsearch} over the acquisition function, with a number of restarts. This derivative free technique can be used to optimize acquisition functions which are not continuous everywhere \cite{DirectSearch1,DirectSearch2}, which is the case for QD with white-box descriptor functions due to the way the acquisition function compares improvement over region-specific elites, making it only piece-wise continuous. The result from each optimisation run is added to a list of \textit{candidate points} and the best point of these is evaluated on our true black-box function.
As the surrogate model of descriptor functions improves --- and therefore posterior uncertainty over region membership reduces --- the boundaries between regions become sharper and approximate the discontinuities of the white-box case above. The acquisition value, at the boundary between regions, can jump dramatically, causing problems for continuous optimizers. As well as using direct search, and in order to avoid issues with getting stuck in local optima, BOP-Elites attempts to identify a large number of diverse starting points in each iteration.

Using Sobol Sampling over $X$ and evaluating each solution on the descriptor surrogate with the posterior mean descriptor, we predict region membership for each point in the sample. The points are evaluated on the EJIE acquisition function - comparing only against the elite indicated in the region prediction - and then ranked. The best performing points predicted to be from unique regions will become the initial restart points for local optimisation along with a number of randomly selected points. The idea here is that we start each optimisation round from promising points well distributed in descriptor space.
The optimisation process is detailed in the pseudo-code for 
Alg.~\ref{pseudocode} and Alg.~\ref{EJIEpseudo}.

\subsection{Slow down of convergence with hard to reach regions}\label{hardtoreach}
Some regions of the descriptor space may be hard to reach because they are very small. It makes sense for an algorithm to spend time seeking out these regions as they contribute both to the diversity of the solution set and the total QD score, but how can we tell the difference between such a region and one that is in fact unreachable?

Our experimental setup, using  MAP-Elites style discrete regions, makes this a very difficult task as the GP surrogate gives non-zero probability to neighbouring descriptor values even when the GP is a good predictive model. As we get close to a discrete boundary in continuous space that probability will increase to nearly 0.5 of the posterior probability, even if the GP is well trained with indicative points. If the neighbouring region is not reachable, then the improvement calculation from expected improvement will be large leading to a huge over-valuation of a point on the boundary of reachable space.

In order to tackle this problem, we introduce a cut-off value to the probability which means that EJIE becomes 

\begin{equation}
EJIE^+(x) = \frac{\sum_{i=1}^{R}\mathbb{\rho}(x \rightarrow R_r)EI_{r_i}(x)}{\sum{\rho}}  \label{EJIE+}
\end{equation}

\begin{equation*}
\rho(x\in r_i) = \begin{cases}
\mathbb{P}(x \rightarrow R_i|D) &\text{if $\mathbb{P}(x\rightarrow R_i|\mathcal{D})>\omega$}\\
0 &\text{otherwise}
\end{cases}
\end{equation*}

The idea here is to ignore regions in our calculations that are too unlikely to be the region of the considered solution candidate. The cutoff value $\omega$ has an initial value of $\frac{1}{(\text{number of regions})}$ which should be easily surpassed by any region we wish to consider when our GP is initialised with a sufficient number of well distributed points. Following \cite{Jones98}, we take that number to be 10 times the number of dimensions of the problem. We expect our confidence in the model to increase as we take more observations so we introduce 

\begin{equation}
\omega = \frac{1}{2}\left(\frac{2}{R}\right)^{\gamma} , \quad \gamma = \sqrt{\frac{10d}{(\alpha - 2\beta +  t)}} \label{cutoff}
\end{equation}
where $\alpha$ is the mis-specification count, i.e., the number of times the algorithm incorrectly attributes more than 50 percent of the acquisition value to a single region which, when evaluated, is a misprediction in descriptor space. $\beta$ is the over-specificity value, the number of times the algorithm is unable to find any expected improvement due to excluding too many regions. $\beta$ therefore increases whenever the optimisation step returns no value in any \textit{candidate points}. This quantity has the desirable property of being $\frac{1}{n}$ when the algorithm starts, due to the initial budget being 10d, and as $t \to \infty , \alpha \to 0$ and therefore $\omega \to \frac{1}{2}$.

\subsection{Dealing with invalid points}\label{invalid}

For many real-world problems, there are regions of the search space that for various reasons do not produce a valid performance or descriptor evaluation. In this work we tried to take the approach that SAIL did, which is to assume that non-convergence is a cheap result and we can essentially treat such events as having no real cost. In this way, if BOP-Elites fails to receive a successful evaluation of a point, it will move on to the next candidate point, or search the acquisition function for an alternative. 

This performed very poorly in our experiments as BOP-Elites searches for points it predicts are high performing and unless it is able to reduce the predicted performance of a point it will continue to suggest it. We thus model validity directly via means of an extra classifying GP or, in our case, a Support Vector Machine (SVM). We chose the latter due to the experimental assumption that NaNs (a Not-a-Number error from the simulator), or invalid points, are cheap or free to evaluate and SVMs work well with a large number of data points, with low computational overhead. 

In addition to classification we implement the Platt scaling method \cite{hearst1998support} to produce a probability estimate for feasibility $\mathbb{P}_{feas}(x)$. We weight our acquisition value by the probability of feasibility. so we now have 

\begin{equation}
EJIE(x)^{++} = EJIE(x)^{+}\mathbb{P}_{valid}(x) \label{EJIE++}
\end{equation}

The modelling of validity becomes active when the first invalid point is found, otherwise BOP-Elites uses $EJIE(.)^{+}$ as standard. 

\subsection{Extending BOP-Elites to work with many regions\label{estfstar}}
When the budget exceeds the number of regions to be optimised, BOP-Elites should work very well in managing exploration and exploitation. However, if the number of regions far exceeds the search budget, BOP-Elites will attempt to focus on high-performing unexplored regions which will maximise the returned points but may give up some exploration for exploitation. SAIL, a surrogate \textit{assisted} algorithm, returns a Prediction Map (PM) as a final output, which is created by running MAP-Elites over the final posterior mean GP surrogate model. The QD-score over the resulting best predicted solutions is used as performance measure. 

BOP-Elites is designed for black box descriptor functions, in which case we additionally model descriptor functions. This introduces the possibility of mis-prediction, i.e., predicting that a point maps to a certain region which, when evaluated on the true descriptor functions, turns out to be incorrect. In the case of mis-prediction a solution contributes no value to the QD score. In the black-box case, the PM may be considered a 'risky' solution as when our surrogate models are not accurate, we expect a severe drop in performance due to mis-predictions. 
BOP-Elites will therefore produce two results. 
\begin{enumerate}
     \item A solution archive made of observed points 
     \item A Prediction Map
 \end{enumerate}
\subsubsection{Prediction Maps} \label{predmaps}
Throughout the optimisation process the performance of the algorithm is reliant on the quality of the surrogate models. Analysing the quality of the underlying models can be useful in seeing how the algorithm is performing and one way of doing such is through a prediction map (PM). A PM uses the surrogate models to predict elites for all regions, including those we have not visited. PM's are created by performing MAP-Elites over the surrogate models where the value function is chosen to find the best points without exploration. In SAIL's original formulation, as the descriptor values are known, the performance measure used to generate the PM is simply the posterior mean of the objective GP and the recommended point for any region $R_i$ should be the posterior mean maximal point for each region:
 \begin{equation}
 \hat{e}_r = \argmax_x \left[\mu|b(x) \rightarrow R_r \right] \label{SPHENPM}.
 \end{equation}
 In the black-box descriptor case we must also predict values for the descriptor. The method implemented in the original SPHEN algorithm is to assume the descriptor GPs are correct and accept the region predicted by the posterior means. In this work we refine this value by acknowledging that neither SPHEN or BOP-Elites actively control exploration of descriptor space, except by virtue of search for high objective performance. In this case we adapt Eqn.\ \ref{SPHENPM} as:
 
\begin{equation}
\hat{e}_r = \argmax_x \left[  \mu|b(x) \rightarrow R_r \right] \mathbb{P}(b(x)\rightarrow R_r ) \label{NEWPM},
\end{equation}

where we multiply the posterior predicted mean value by the probability of x being in the region.
 
\subsection{Upscaling of archives}
One interesting unexplored potential of the PM approach is that while a structured archive imposes a discretisation over the descriptor space, the models learned throughout the optimisation process are continuous. This allows us to use points collected during a coarsely discretised solution archive to build a rich model of the underlying functions and predict solutions for a much more finely discretised solution archive. We call this process 'upscaling'.

A key difference between previous surrogate assisted methods and BOP-Elites is the way they balance exploration and exploitation. SPHEN uses a tunable parameter $\beta$ to weight towards variance reduction, which leads to early model-building behaviour\footnote{In this paper we present results where the $\beta$ has been tuned, but incorrectly specified $\beta$ can lead to poor results.}.  By contrast BOP-Elites is improvement greedy which means it fills regions starting from the highest performing regions first. As BOP-Elites takes many observations in a concentrated area the underlying models can suffer until BOP-Elites gradually moves away from the higher performing regions. As SPHEN reduces posterior variance early, it leads to a good early model.

A good early model is desirable for BOP-Elites, as the quality of predictions and therefore performance of algorithm depend on the model. In this case we propose to utilise the region filling behaviour of BOP-Elites as a tool.

\subsubsection{Initial Upscaling}
As BOP-Elites is excellent at filling regions with high performing points, we first give it the task of filling a coarsely discretised solution archive, say a two dimensional descriptor space of $5^2$ regions, $[5,5]$. BOP-Elites will fill these regions, creating diversity over the descriptor space due to how spread out the regions are. Once BOP-Elites has either filled the regions or has spent a small budget, twice the number of regions, we increase the discretisation to the original resolution, leveraging the now improved model. This initial phase of BOP-elites improves the early prediction models and improves the surrogate models aiding the optimisation process throughout.

\subsubsection{Upscaling to higher dimensions}
As the number of points a Gaussian Process model can feasibly handle is limited due to computational complexity and numerical issues, QDO problems with higher resolution solution archives are infeasible for surrogate based methods to process naively. Whilst possible solutions may include sparse GPs, inducing points etc., in this work we take a simple approach to make predictions for higher resolution solution archives from optimised lower resolution solution processes. The idea here is to take the full observation set from a $25\times25$ solution (both the elite set and the observation set used in the optimisation) and make continuous predictions from the GP to fill a higher resolution grid. We present the results of these efforts in the following sections.

\colorlet{transgray}{gray!10}
\sethlcolor{transgray}

\begin{algorithm}
\begin{algorithmic}[1]
\Statex \textbf{Initialization:}
\State Sample initial dataset $\mathcal{D}_0 = \{ X_{0},Y_{0},\textbf{b}_{0} \} $ 
\State Elite points in each region enter solution archive, $S $ 
\State Build Objective and Descriptor surrogate models $\hat{f}(.) , \hat{\textbf{b}}(.)$
\Statex
\Statex \textbf{Optimization Loop:}
\While {$n < N_{max}$}
\For{points in $X_d$}
\State compute all region improvements (Eqn.~\ref{nicheimp})
\State compute all region probabilities   (Eqn.~\ref{nicheprob})
\State compute $EJIE^+(x)$ (Eqn.~\ref{EJIE})
\EndFor
\State sample $x_n$ at $argmax(EJIE(x))$
\If{$b(x_n) \in r_i$ and $f(x_n) > f(\hat{e}_{r_i})$}
\State $\hat{e}_{r_i} = x_n$
\EndIf
\State update posterior models
\State $n = n + 1$
\EndWhile
\State \textbf{Return:} Elites
\caption{BOP-Elites Pseudo-code} \label{pseudocode}
\end{algorithmic}
\end{algorithm}



 

\begin{algorithm}
\begin{algorithmic}[1]
\State \textbf{Inputs}: $x , \hat{f}(.) , \hat{b}(.), \omega $
\State $ \hat{b}(x) \rightarrow \hat{r}_x$ \Comment{predict region membership from posterior mean}

\For{$r$ in $R$}
    \State $EI_r \left( x,\hat{f}_r(.), \hat{e}_r \right) $ \Comment{by Eqn. \ref{EIR}}
    \State $\P_r = \P(x \in r|D)$ \Comment{by Eqn. \ref{probr}}
    \If{$\P_r < \omega$ }
        \State $\P_r$ = 0 \Comment{remove region if probability $<$ $\omega$ as Section \ref{hardtoreach}}
    \EndIf
\EndFor
\State Normalise probabilities
\State EJIE$ = \sum_r\P_r EI_r$
\State \Return EJIE
\caption{EJIE+ Acquisition Function}\label{EJIEpseudo}
\end{algorithmic}
\end{algorithm}

\section{Evaluation Methodology\label{sec:Methods}}
We compare BOP-Elites against state-of-the-art QD algorithms on a number of benchmarks designed to simulate problem settings with interesting properties, including invalid solutions, different dimensionality, and either 1 or 2 descriptor functions. We measure performance by the QD score of the solution archive and when evaluating PMs we observe the true values of the recommended points and calculate the QD score, attributing a default low valuation for contributions from mis-specified points in the black-box case.    



\subsection{Benchmark algorithms\label{problem_maker}}
BOP-Elites is designed for problems where descriptor functions are black-box and assumed to be expensive. We compare against 3 baselines: Sobol sampling, MAP-Elites and the state of the art surrogate assisted method which is SAIL for the white-box descriptor case and SPHEN when the descriptors are black-box functions.

Sobol sampling \cite{Sobol1976,SobolScrambled} is a quasi-random number generation technique designed to sequentially generate a set of points that converge to uniformity over the unit hypercube as the number of samples $t$ approaches infinity. While Sobol sampling is not a Quality-Diversity algorithm per se, retaining the best points found in each region can be employed to construct a valid solution archive.

We use MAP-Elites \cite{cully2015robots} with uniformly sub-partitioned descriptor dimensions. In our experiments we use Gaussian mutation with mutation step size $\sigma=0.1$ and 50 children per generation, no crossover and the population is the set of elites in the archive.

For problems with white-box descriptor functions, we compare against the SAIL algorithm \cite{gaier2018data} which builds a surrogate model of the objective with exactly the same specification as BOP-Elites' GPs using the BOTORCH \cite{balandat2020botorch} framework --- see the next section - and runs MAP-Elites over those GP's with the same setup as standalone MAP-Elites, as above, except that it uses the UCB function as the objective. This function has a user-defined parameter $\beta$, and we have optimised it at $\beta = 3.7$ which performed well across all region numbers and best, on average, across domains.

SPHEN is an implementation of SAIL, with all the same settings, that additionally models the descriptor values with a surrogate model as described in Section~\ref{SPHEN}. 


\subsection{Experimental configuration}
BOP-Elites uses a noise-free GP model implemented in BOTORCH using a Matern 5/2 kernel with standardised objectives and normalised search inputs. We use the pattern search optimiser by PYMOO in the search step with default settings and a maximum budget of 100 generations of 1500 evaluations. BOP-Elites has one parameter which may benefit from problem specific tuning; the number of restarts per iteration in the optimiser. Rather than fine-tuning this, we chose a naive setting of 10 restarts for the $10\times10$,$25\times25$ and $50\times50$ solution archives, increasing this budget is almost certain to further improve the performance of BOP-Elites.
We compare BOP-Elites with a small budget of observations against MAP-Elites at different observation budgets. As MAP-Elites generates 50 children per generation each 1000 observations equates to 20 generations. 

Our results compare BOP-Elites performance with both the observed solution archive, the final set at the end of the run, and the value of the predicted map, the value of the predicted points made after optimisation. We chose to limit the number of observations for BOP-Elites to 1000 for the $10\times10$ and 1250 for the $25\times25$ and $50\times50$ archives due to the computational time of making predictions from the Gaussian process at high observation count. 1250 is double the number of regions in the 25 by 25 structured archive and was sufficient to find a very strong solution. For the $50\times50$ archive, BOP-Elites does not have sufficient budget to create a good solution but we include these results to show a failure case which can be resolved by upscaling.


\subsection{Benchmark functions}
\subsubsection{Synthetic Gaussian Process problems}
We generate a single output, 10 dimensional input, Gaussian Process model with a predefined lengthscale and take one realisation of the posterior Gaussian Process as a deterministic objective function for testing. We are able to generate both objective and one or more descriptor functions this way, combining them to create multidimensional descriptor spaces.

The benefit of generating functions from GPs is that we can create an infinity of unique functions quickly and efficiently without concern for model mismatch with the GP surrogate model. Lengthscales are chosen uniformly at random from 0.1 to 1 and algorithms are compared on functions generated from the same seeds.

\subsubsection{Mishras Bird function}
Mishras Bird function is an optimisation benchmark \cite{MishraBirdFn} with the following analytical form : 
\begin{align*}
f(x_1,x_2) &= A + B + C, \\
A &= \sin{x_2}\exp(1-\cos{x_1}^2), \\
B &= \cos{x_1}\exp(1-\sin{x_2})^2, \\
C &= (x_1-x_2)^2, \\
\text{with } & -10\le x_1 \le 0, -6.5 \le x_2 \le 0.
\end{align*}

We define 2 simple descriptor functions with Mishras bird function:
\[B_1(x) = -x_1 , B_2(x) = -x_2\]

Whilst being low dimensional, the problem exhibits an interesting property as there is a small region that vastly outperforms the rest of the domain in terms of function value.

\subsubsection{Robot Planar Arm}
This well-known QD benchmark problem \cite{cully2015robots} is a simulation of a planar robot arm with an $n= 4$ dimensional unit hypercube as the input and objective and descriptor values as follows:

\begin{eqnarray}
y &=& 1- \sqrt{\frac{1}{n}\sum_{i=1}^{n} (x_i-\bar{x})^2} \\
B_1 &=& \frac{1}{2n}\sum_{i=1}^n \sin\left(\sum_{j=1}^{i}{( 2\pi x_j- \pi)}\right) +0.5 \\
B_2 &=&\frac{1}{2n}\sum_{i=1}^n \cos\left( \sum_{j=1}^{i}{(  2\pi x_j- \pi)}\right) +0.5
\end{eqnarray}

where $\bar{x}$ is the mean of $x$, $\bar{x} = \frac{1}{n}\sum_{i=1}^{i=n}x_i$, the behaviour space is known to be $[B_1,B_2] \in [0,1]^2$ and the objective function has values in the range $y \in [0,1]$

This problem simulates the ability of a robot arm to reach feasible regions in a $2d$ plane while minimising the standard deviation of the position of the individual joints.

\subsubsection{Rosenbrock 6d}
The Rosenbrock function \cite{rosenbrock} is defined for even numbered dimensions, here we use 6 dimensions. 
\[f(x) = \sum_{i=1,3,5} 100\left[ (2x_i-2x_{i+1}^2)^2 + (1-2x_{i+1}^2) \right] 
\]
\[0 \le x_i \le 1 , \forall i.\]

We used 2 feature dimensions with the form 
\[B_1(x) = \frac{1}{2}(x_1+x_2) , B_2(x) = (x_3 -1)^2\]

Rosenbrock is a relatively smooth function but in 6 dimensions.

\subsubsection{PARSEC: 2D Aerofoil optimisation problem}
A 2D aerofoil design was modelled with the Xfoil simulation package \cite{xfoil} and the objective and features were derived from both the physical characteristics of the aerofoil and the simulation results. This is a well known and well studied problem in engineering and we used the 10 dimensional PARSEC representation of the parameters \cite{sobieczky1999parametric}. This problem provides an informative baseline as it forms a key evaluation from the previous SAIL papers. The objective is calculated with penalty values $\rho$ as:
\[ objective(x) = \left( \mu_{drag(x)} + K\mu_{drag(x)}\right) \times \rho_{lift}(x) \times \rho_{area} (x)\]
and the descriptors are values taken from the parameterised representation called $X_{up}$ and $Z_{up}$. For specifics on these values we guide readers to the SAIL literature. 

In this work, BOP-Elites is configured similarly to the prior SAIL research. The approach utilizes the EI calculated exclusively from the GP modeling drag, with a lift penalty also modeled with a GP being applied to the EI value post-calculation. This unique formulation uses the lift-associated penalty more as a classifier, rather than allowing it to influence the performance prediction used in the EI computation. As a result, the exploration is not driven by the variance in the lift.
The PARSEC domain features invalid regions, for which the adaptations in Section~\ref{invalid} were implemented.

\section{Empirical results\label{sec:results}}

We present tables of summary statistics from 100 experiments for each algorithm. In the top row of a results table we compare the performance of the evaluated solutions in the archive of algorithms after 1000 function evaluations in the 10x10 case or 1250 function evaluations for the 25x25 and 50x50 case. Besides the evaluated solution archives we generate PMs, which are generated as in Section \ref{predmaps} and measured as the true value of the points suggested as in Eqn. \ref{TEP}. We compare these results alongside MAP-Elites and SOBOL sampling using a 50k evaluation budget, which appear in a separate row to clarify that a much higher budget was used. Highlighted cells indicate the best value found among all solutions. Yellow cells indicate the best mean results, statistically tested above a 99\% confidence level with a t-test. In cases where the same algorithm generates both an observed and a predicted result, and neither of these results significantly outperforms the other, we only denote the result that exhibits the highest mean performance. 

\subsection{Black box descriptors}
\subsubsection{Synthetic Gaussian Process functions - blackbox descriptor functions}
\begin{figure}
    \centering
   \includegraphics[width=0.48\textwidth]{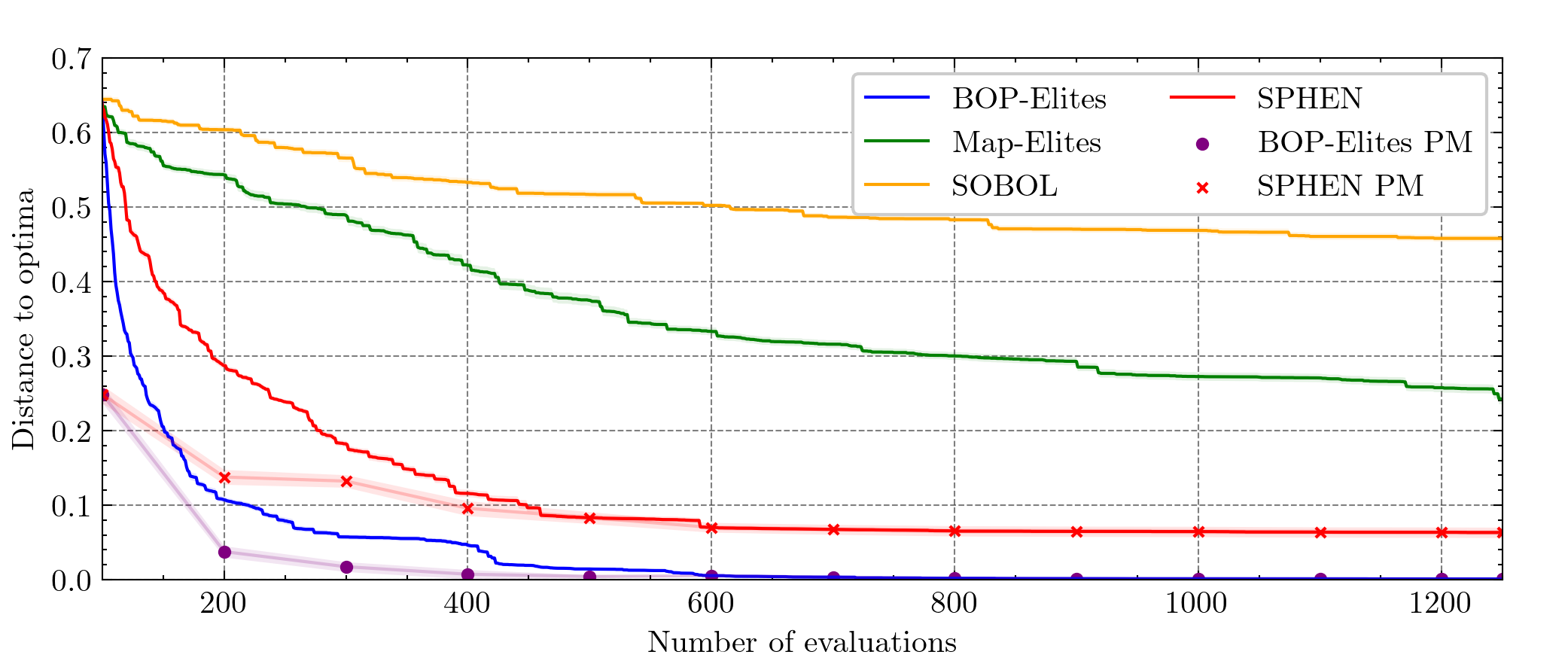}
    \caption{Convergence of BOP-Elites in comparison to Map-Elites, SOBOL and SPHEN on the synthetic GP problem with a 1d 25 region archive. We normalise the performance of the algorithms for each seed before taking the mean over 100 experiments.\label{fig:synth_convergence}}
    
\end{figure}

Table~\ref{tab:syntheticBlack} contains results from experiments on synthetic GP functions with 1 generated descriptor function. As each random seed produces a unique objective and feature function, we normalise the performance (against the best score found by all algorithms) for each seed and take the average. BOP-Elites significantly outperforms all other algorithms even at much higher function evaluations. SPHEN appears to converge to a sub-optimal score and this is because it fails to fill all available regions, particularly in the 10 region case. Figure \ref{fig:synth_convergence} shows both SPHEN and BOP-Elites increase in QD-score value quickly, with BOP-Elites improving quickest and converging at an optimal value. Synthetic GP benchmarks help to eliminate the problem of model mis-specification for surrogate based methods and this showcases the potential of BOP-Elites in optimal conditions. While SPHEN is also a surrogate based method it struggles to achieve optimal performance and examination of the returned values show that it sometimes struggles to reach all the regions.
\begin{table}
\caption{Final QD score on the 10d synthetic benchmark with black box descriptor functions. Algorithms on the top use 1000/1250 function evaluations for respective archive sizes. Algorithms at the bottom use 50k function evaluations. \label{tab:syntheticBlack}}
\begin{center}
\begin{tabular}{r@{\qquad}rr@{\qquad}rr}
\hline
QD score:  & \multicolumn{2}{c}{1d. 10} & \multicolumn{2}{c}{1d. 25} \\
\cline{2-5}
 & Mean & S.E. & Mean & S.E.   \\
\hline
BOP-Elites &  \cellcolor{gray!25} 1.00 & 0.000  &  \cellcolor{gray!25}1.00& 0.0005   \\
SPHEN &  0.73 & 0.007  &  0.88& 0.0052\\
BOP-Elites-PM  & \cellcolor{gray!25}1.00 & 0.000  & 0.99& 0.0023  \\
SPHEN-PM &  0.76 & 0.008  &  0.88& 0.0052\\
MAP-Elites  & 0.63 & 0.011  & 0.72& 0.0087 \\
\hline
MAP-Elites 50K  & 0.93 & 0.005  & 0.97& 0.0091  \\
Sobol 50K  & 0.81 & 0.009  &0.76 & 0.0103 \\
\hline
\end{tabular}
\end{center}
\end{table}






\subsubsection{Mishras Bird function - blackbox descriptor functions}
\begin{figure}
    \centering
    \includegraphics[width=0.48\textwidth]{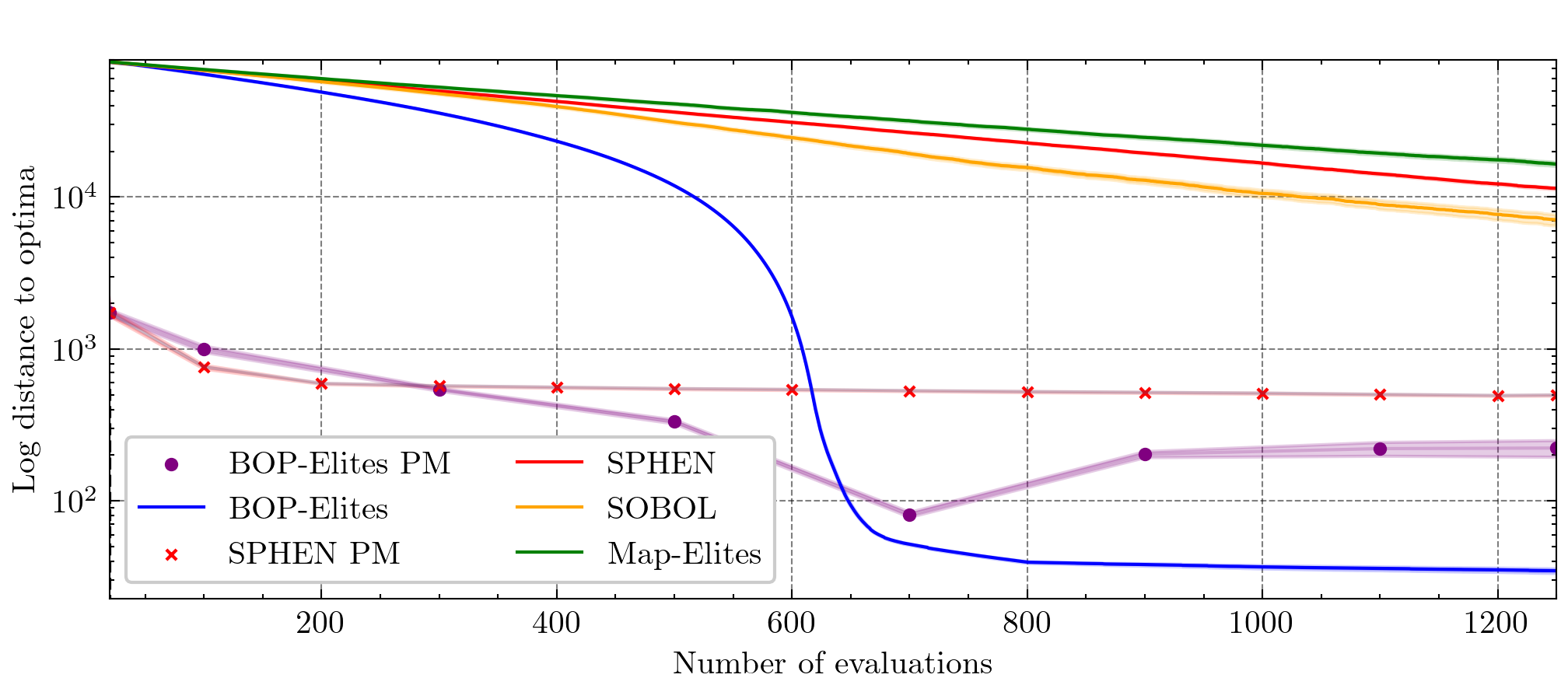}
    \caption{Convergence of BOP-Elites in comparison to Map-Elites and SOBOL for the Mishra Bird Function with 2 descriptor dimensions and a $[25 \times 25]$ solution archive. BOP-Elites rapidly finds a near optimal solution. \label{fig:Mishra_convergence}}
\end{figure}

In Table~\ref{tab:mishraBlack} and Fig.~\ref{fig:Mishra_convergence}, we see strong early performance of BOP-Elites for both prediction map solutions. This is exactly as we would expect with a smooth and well behaved combination of objective/descriptor functions. BOP-Elites clearly outperforms other methods in the final stages, and the better points found by BOP-Elites enable an improvement in the PMs. We notice two things in the performance of the BOP-Elites PM: The error is lowest at 700 function evaluations before increasing again, we also notice an increase in the Standard Error of the PM performance. While the error is relatively small, this unusual behaviour is almost certainly due to the way BOP-Elites is focusing on visiting high performing points and not on reducing the uncertainty of the model, in particular in the descriptor predictions. 

Another intriguing observation is the impressive early performance of SOBOL-sampling in this domain, surpassing that of SPHEN. This outcome can be primarily attributed to the smooth, flat behavior across the function, with value distribution being even throughout the domain. As a result, effective performance can be achieved simply by covering the input space, an aspect that SOBOL excels at in two dimensions where the input also serves as the descriptor.
\begin{table}
\caption{Final QD score on the Mishra 2d benchmark with black box descriptor functions. \label{tab:mishraBlack}} 
\begin{center}
\begin{tabular}{r@{\qquad}rr@{\qquad}rr}
\hline
QDS: Mishras Bird & \multicolumn{2}{c}{10x10} & \multicolumn{2}{c}{25x25}  \\
\cline{2-5}
 & Mean & S.E. & Mean & S.E.  \\
\hline
BOP-Elites & 13776.18 & 0.45  &\cellcolor{gray!25} 79973.29 & 1.186 \\
SPHEN   & 13762.26 & 3.87  & 68587.21 & 237.02 \\  
BOP-Elites-PM  & \cellcolor{gray!25}13776.55 & 0.28  & 79785.12 & 27.8 \\
SPHEN-PM&13744.15 & 1.44  & 79552.95 & 0.81 \\
MAP-Elites  & 11136.70 & 89.61  &25473.58 &196.45  \\
\hline
MAP-Elites 50K & 13666.59 & 1.63  & 79033.73 & 5.80  \\
Sobol 50K & 13710.03 & 1.08  & 79536.39 & 2.49  \\
\hline
\end{tabular}
\end{center}
\end{table}






\subsubsection{Robotarm - blackbox descriptor functions}
\begin{figure}
    \centering
    \includegraphics[width=8cm]{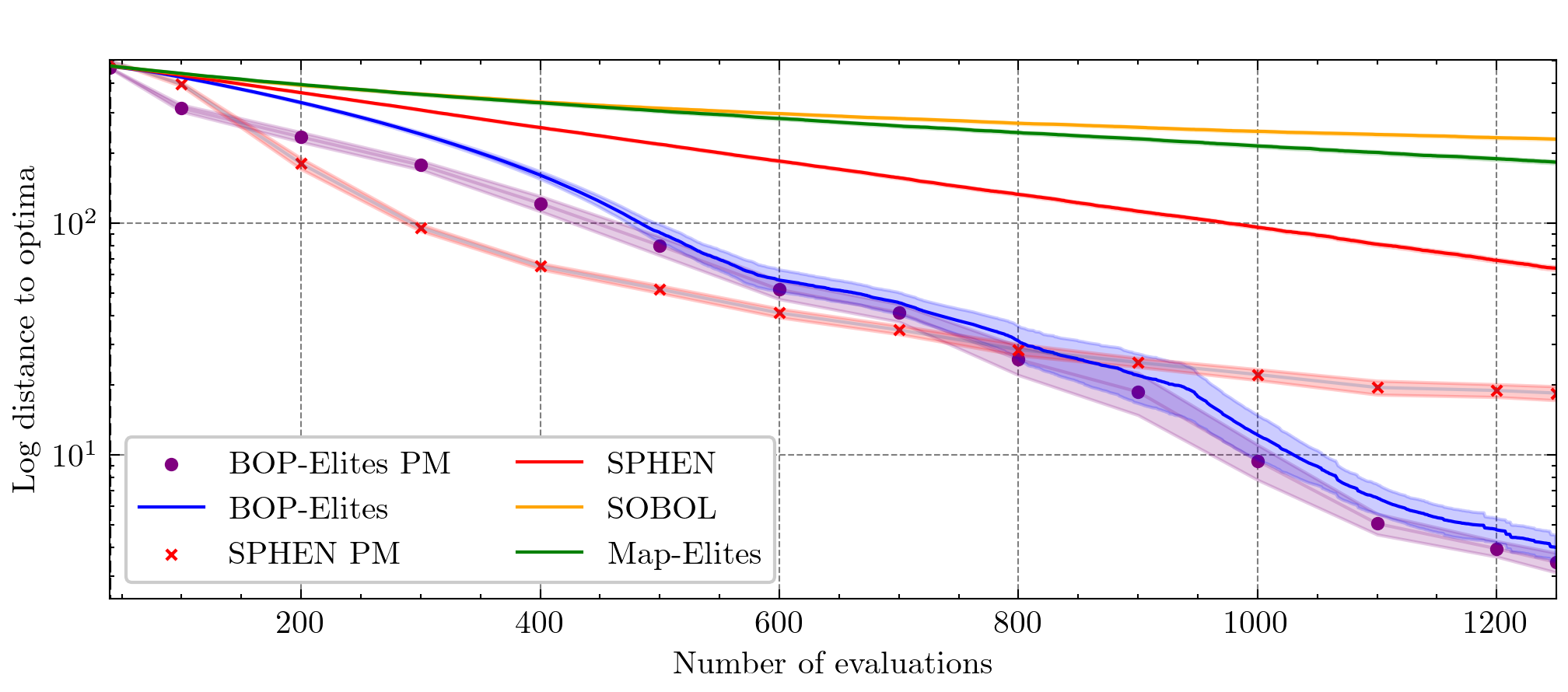}
    \caption{Convergence on the Robot Planar Arm function with 2d black-box descriptor space over a $[25 \times 25]$ solution archive. A comparison between BOP-Elites and Map-Elites, SOBOL and SPHEN. BOP-Elites and SPHEN model both objective and descriptor values. Lines represent mean performance with the standard error around the mean.    \label{fig:robotarm_convergence_BB}}
\end{figure}

Table~\ref{tab:robotBlack} summarises results on the robotarm benchmark problem. BOP-Elites returns the best solution archive for the 10x10 and 25x25 solution archive and much better than even MAP-Elites with 50k evaluations. Figure~\ref{fig:robotarm_convergence_BB} shows BOP-Elites has excellent early performance that slows between 600 and 900 evaluations before accelerating. This is likely where BOP-Elites is refining the $\omega$ value used to focus its search and filter out unreachable regions. The excellent end performance indicates that BOP-Elites is also a good choice for QD problems with unreachable regions. SPHEN doesn't require any adaptations to handle the unreachable regions in this domain and shows very good early performance, but falls behind BOP-Elites towards the end of the budget.

\begin{table}
\caption{Final QD score on the robotarm benchmark with black box descriptor functions. .
\label{tab:robotBlack}}
\begin{center}
\begin{tabular}{r@{\qquad}rr@{\qquad}rr}
  \hline
  QD scores: robotarm  & \multicolumn{2}{c}{10x10} & \multicolumn{2}{c}{25x25} \\
\cline{2-5}
  & Mean & S.E. & Mean & S.E.    \\
  \hline
  BOP-Elites & \cellcolor{gray!25}85.14 & 0.00  & 500.12 & 0.1 1\\
  SPHEN &  84.48 & 0.68  & 464.41 & 8.485   \\
  BOP-Elites-PM   &84.91 & 0.02  &\cellcolor{gray!25} 502.30 & 0.14 \\
  SPHEN-PM   &83.49 & 0.17  &451.07 & 0.26 \\
  MAP-Elites  & 50.27 & 0.69  & 140.10 & 1.31  \\
  \hline
  MAP-Elites 50K & 84.15 & 0.02  & 493.15 & 0.40  \\
  Sobol 50K & 81.21 & 0.08  &267.07 & 1.15  \\
  \hline
\end{tabular}
\end{center}
\end{table}

\begin{figure*}[!t]
    \centering
    \subfloat{\includegraphics[width=0.3\linewidth]{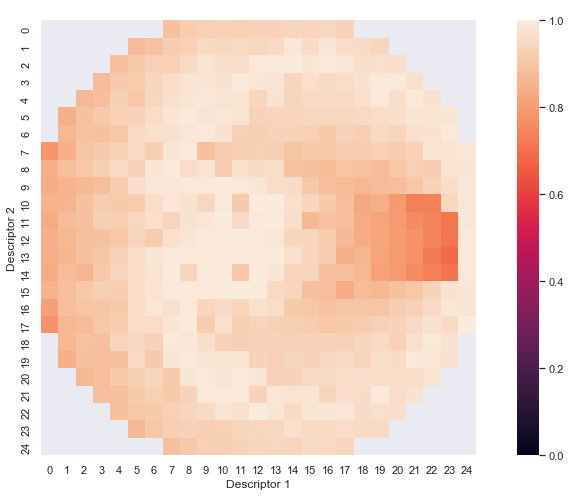}} 
    \hfill
    \subfloat{\includegraphics[width=0.3\linewidth]{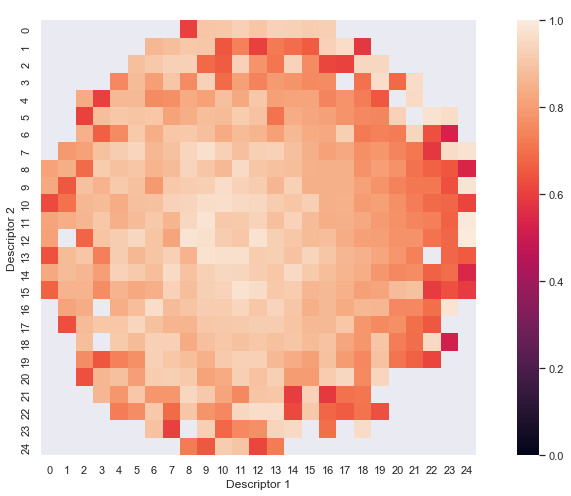}} 
    \hfill
    \subfloat{\includegraphics[width=0.3\linewidth]{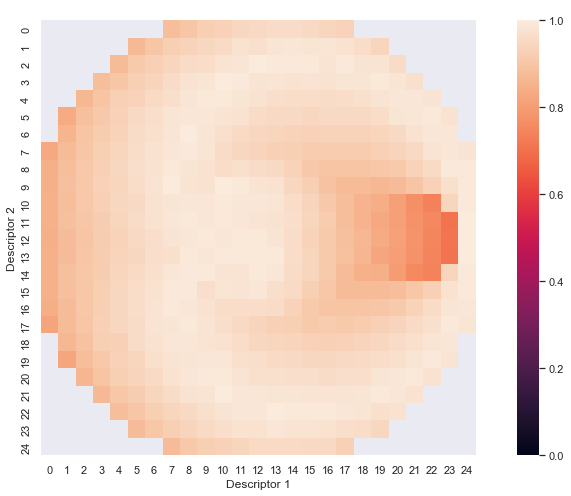}}
    \caption{These plots show an indicative result for a single solution archive of 25x25 : (left) BOP-Elites finds a near-optimal solution and explores all regions after 1250 observations. (centre) After 10,000 observations, MAP-Elites is only starting to fill all regions and find high-quality solutions. (right) MAP-Elites, with 200k total observations, achieves a near-optimal solution set, surpassing BOP-Elites only after $\sim{100k}$ points.}
    \label{fig:robotarm25}
\end{figure*}

\subsubsection{Rosenbrock 6d function - blackbox descriptor functions}
\begin{figure}
    \centering
    \includegraphics[width=0.45\textwidth]{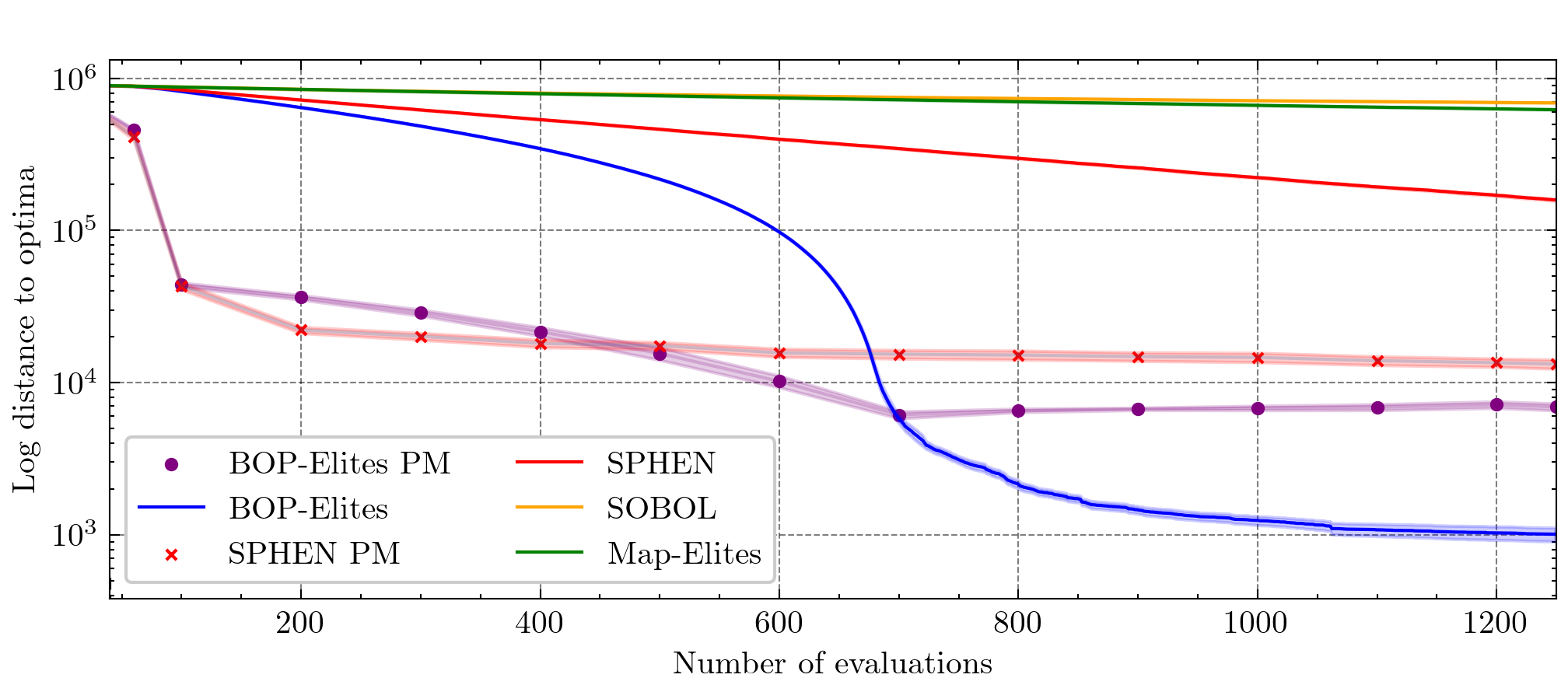}
    \caption{Convergence of BOP-Elites in comparison to Map-Elites and SOBOL on the Rosenbrock6d problem with a 2d black-box descriptor space and a $[25\times25]$ solution archive   \label{fig:rosenbrock_convergence}}
\end{figure}

Table~\ref{tab:rosenbrocktable}  and Fig.\ref{fig:rosenbrock_convergence} shows a comparison of BOP-Elites, SPHEN, MAP-Elites and Sobol sampling for the Rosenbrock function with 6 input dimensions. We again see dominant performance for BOP-Elites and corresponding strong performance for the BOP-Elites PM, with an interesting dip in the PM error at around the 700 observations mark. At this point in the algorithm BOP-Elites changes its behaviours from exploration (filling regions) to exploitation (improving existing regions). As regions, by definition, represent diversity, the PM is now not improving as rapidly as the points found are no longer as useful in improving the surrogate model. 

Additionally, the performance of the BOP-Elites Prediction Map now under-performs, on average, the observed solution map towards the end of the budget. This appears to be due to uncertainty over region membership at the boundaries leading to occasionally suggesting a point that ends up being in a region other than that predicted, leading to performance loss. 

\begin{table}
\caption{Final QD score on the Rosenbrock 6d function with black box descriptor functions.  }\label{tab:rosenbrocktable}
\begin{center}
\begin{tabular}{r@{\qquad}rr@{\qquad}rr}
\hline
QDS: Rosenbrock6d & \multicolumn{2}{c}{10x10} & \multicolumn{2}{c}{25x25}  \\
\cline{2-5}
& Mean & S.E. & Mean & S.E.   \\
\hline
BOP-Elites & \cellcolor{gray!25}152837.88 & 0.04  & \cellcolor{gray!25}903547.61 & 2.27 \\
SPHEN   & 151488.48 & 518.90  &822086.93 & 8011.94\\       
BOP-Elites-PM  &152340.82 & 48.27  &903327.12 & 4.20\\
SPHEN-PM  &152450.18 & 48.27  &896067.23 & 3.57 \\
MAP-Elites & 99826.52 & 970.44  &254849.92 &2708.37  \\
\hline
MAP-Elites 50K & 151674.80 & 18.31  & 899302.84 & 233.07  \\
Sobol 50K & 106240.71 & 188.91  & 543038.19 & 344.11  \\
\hline
\end{tabular}
\end{center}
\end{table}

\subsection{White Box Descriptors}

For problems with White-Box descriptors, we compare against the SAIL algorithm rather than SPHEN. BOP-Elites is given the same budget as in the Black-Box case but for the 10x10 and 25x25 problems it often converges to a good solution early. PMs in this case are created with known descriptor functions.

\subsubsection{Mishra Bird function - white-box descriptor functions}
Table~\ref{tab:mishrawb} shows summary statistics for the Mishra Bird function with white-Box descriptors. BOP-Elites finds a near optimal solution (Comparably closer to the optimal than in the Black-Box case). Perhaps surprisingly SAIL does not appear to benefit much from knowledge of the descriptor functions, meaning the inability to keep up with BOP-Elites is likely down to the randomness in SAILs point selection mechanism.

\begin{figure}
    \centering
    \includegraphics[width=0.48\textwidth]{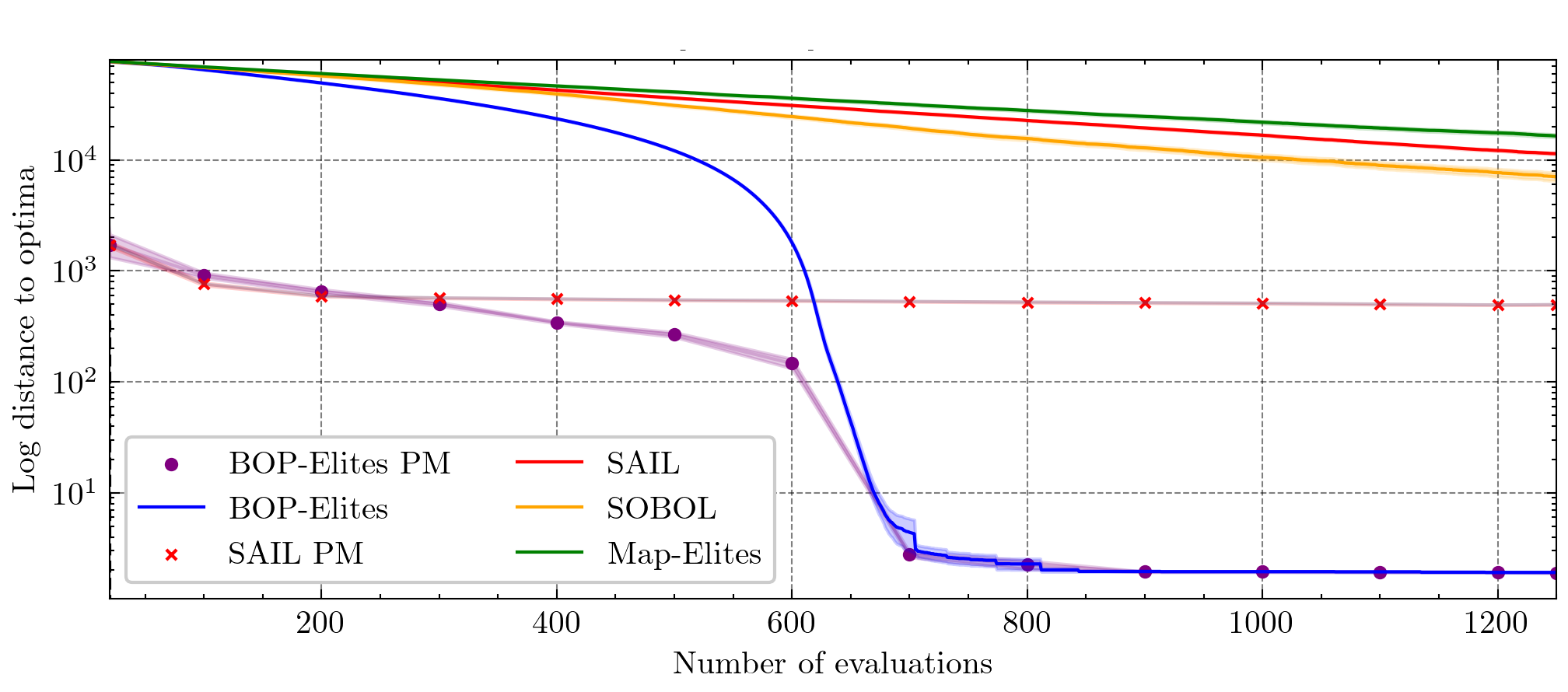}
    \caption{Convergence of BOP-Elites in comparison to Map-Elites and SOBOL on the Mishra Bird function with 2 white-box descriptor functions over a $[25 \times 25]$ solution archive
    \label{fig:mishra_convergencewb}}
\end{figure}

BOP-Elites however produces a very different quality of PM which keeps step with the observed solution set as seen in Fig \ref{fig:mishra_convergencewb}. This implies the drop in the black-box performance is due to errors in the descriptor models.

\begin{table}
\begin{center}
\caption{Final QD score on the Mishra 2d function with White box descriptor functions.  }\label{tab:mishrawb}
\begin{tabular}{r@{\qquad}rr@{\qquad}rr}
\hline
QDs: mishra problem & \multicolumn{2}{c}{10x10} & \multicolumn{2}{c}{25x25}  \\
\cline{2-5}
& Mean & S.E.. & Mean & S.E..   \\
\hline
BOP-Elites & \cellcolor{gray!25}13778.69  & 0.04   &\cellcolor{gray!25}80005.96 & 00.01\\
SAIL & 13760.23 & 2.45 & 68587.44  & 51.98   \\
BOP-Elites-PM   & 13778.68 &  0.01 &  79995.04 & 0.17 \\
SAIL PM  & 13752.21 & 0.63  & 79561.09 & 1.21  \\
MAP-Elites  & 11136.70 & 89.61  &25473.58 &196.45  \\
\hline
MAP-Elites 50K  &13666.59& 1.63& 79033.73 &05.80  \\
Sobol 50K &13710.03 &1.08 &79536.39& 02.49 \\
\hline
\end{tabular}
\end{center}
\end{table}

\subsubsection{Robotarm - white-box descriptor functions}
In the white-box descriptor case, BOP-Elites does not need to actively model the reachability of regions and this explains why we see such good performance. We do not see the slow down in performance as in the black-box case. As Table \ref{tab:robotarmwb} shows, in this domain the PMs perform slightly better than their fully evaluated counterparts and BOP-Elites performs best in all scenarios. 

\begin{table}
\begin{center}
\caption{Final QD score on the Robot Arm function with White box descriptor functions. \label{tab:robotarmwb} }
\begin{tabular}{r@{\qquad}rr@{\qquad}rr}
\hline
QDs: robotarm problem & \multicolumn{2}{c}{10x10} & \multicolumn{2}{c}{25x25}  \\
\cline{2-5}
& Mean & S.E. & Mean & S.E.  \\ 
\hline
BOP-Elites &\cellcolor{gray!25} 85.17 & 0.001  &504.30 & 0.020  \\
SAIL & 83.98 & 0.300  & 484.13 & 2.510  \\
BOP-Elites-PM  &\cellcolor{gray!25} 85.17  & 0.001  &\cellcolor{gray!25} 505.10& 0.005 \\
SAIL - PM & 84.42  & 0.216 & 499.39 & 0.120\\
MAP-Elites  & 50.27 & 0.690  & 140.10 & 1.310  \\
\hline
MAP-Elites 50K  & 84.15 &0.020 &493.15& 0.400 \\
Sobol 50K  & 81.21 &0.080 &267.07 &1.150 \\
\hline
\end{tabular}
\end{center}
\end{table}

\begin{figure}
    \centering
    \includegraphics[width=0.48\textwidth]{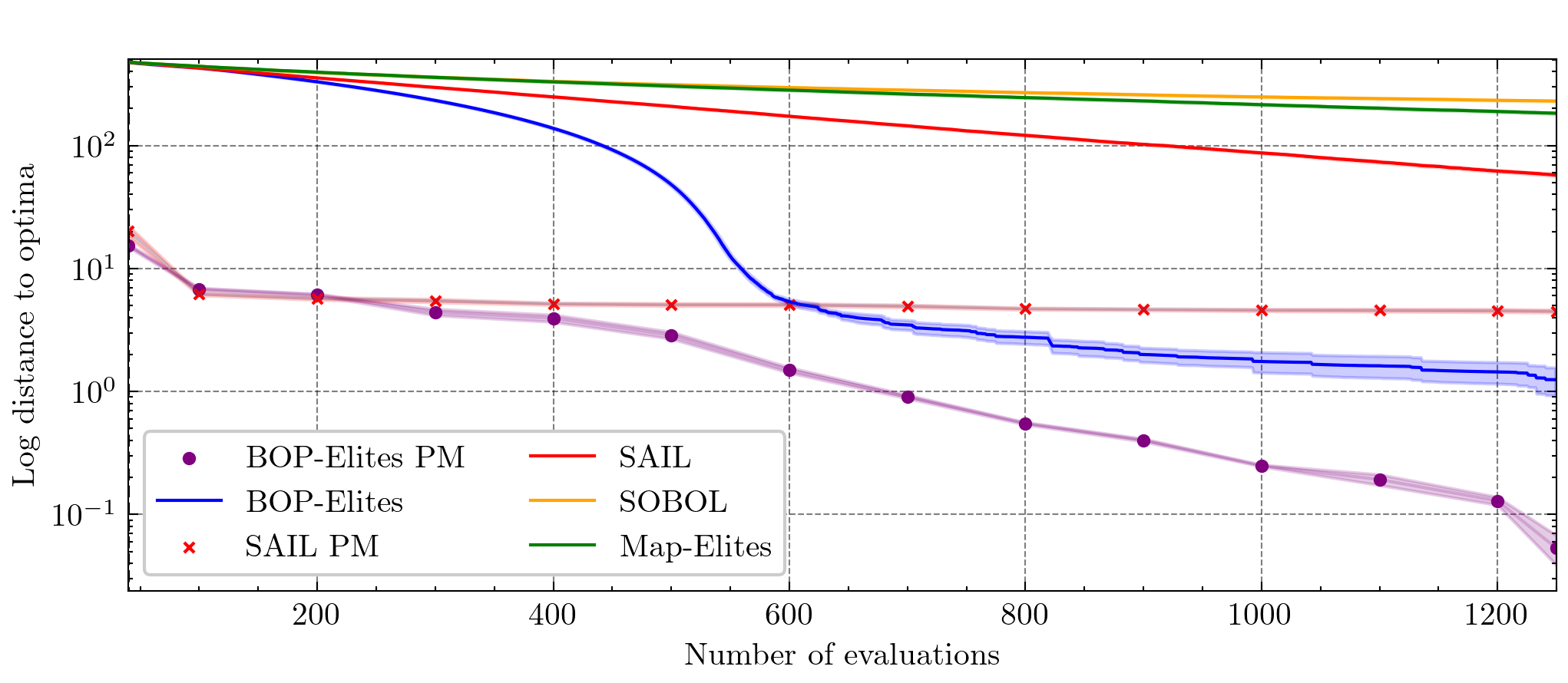}
    \caption{Convergence of BOP-Elites in comparison to Map-Elites and SOBOL on the 4d Robot Planar arm problem with 2 whitebox descriptor functions and a $[25\times25]$ solution archive
    \label{fig:robotarm_convergence}}
\end{figure}

\subsubsection{PARSEC - white-box descriptor functions}
For the PARSEC domain which contains invalid solutions, Table~\ref{tab:parsec} shows that BOP-Elites hugely outperforms the other algorithms. This is unsurprising as BOP-Elites is the only method capable of actively modelling the infeasible regions. While we assume there is no penalty for attempting to evaluate an invalid point, we note that BOP-Elites attempts to make an average of 83 evaluations of invalid points, compared with 373 from the SAIL algorithm from a run of 1250 points in the 25x25 case. As this domain closely models a real-world optimisation task, these results provide an important confirmation of BOP-Elites' efficacy and ability to extend to problems with infeasibility restrictions.

\begin{table}
\begin{center}
\caption{Final QD score on the Parsec domain with White box descriptor functions. \label{tab:parsec}}
\begin{tabular}{r@{\qquad}rr@{\qquad}rr}
  \hline
  QD scores: PARSEC & \multicolumn{2}{c}{10x10} & \multicolumn{2}{c}{25x25} \\
  \cline{2-5}
  & Mean & S.E.  & Mean & S.E. \\
  \hline
  BOP-Elites         &371.01 & 0.21 & \cellcolor{gray!25}2207.42 & 1.17 \\
  SAIL               & 361.21 & 1.63 & 1832.77  & 12.31   \\
      SAIL - PM        & 369.52 & 0.51 & 1875.06 & 1.79  \\
          BOP-Elites - PM  & \cellcolor{gray!25}371.51 & 0.17 & 2198.79 & 1.30  \\
    
  MAP-Elites        & 205.77 & 1.21 & 293.12  & 37.27  \\
  \hline
  MAP-Elites 50k         & 333.54 & 2.21 & 2059.32  & 47.27 \\
  Sobol 50k             &  241.93 & 5.21 & 1341.07 & 34.21 \\
\hline
\end{tabular}
\end{center}
\end{table}

\begin{figure}
    \centering
    \includegraphics[width = \linewidth]{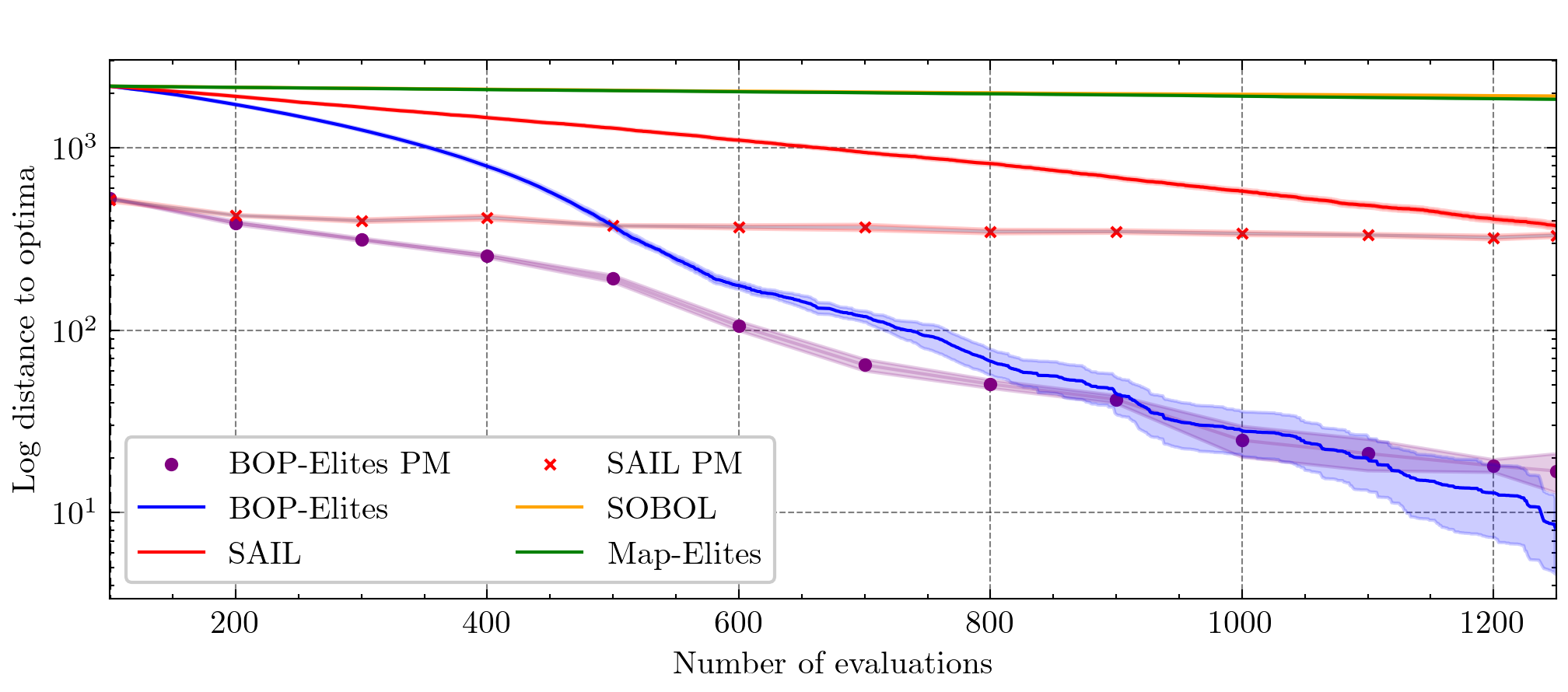}
    \caption{Convergence of QDscores for BOP-Elites in comparison to Map-Elites, SOBOL and SAIL on the PARSEC problem with 2 white-box descriptor functions, a $[25\times25]$ solution archive. }
    \label{fig:PARSECconvergence}
\end{figure}

\subsection{Upscaling Results}
The following tables present results when performing 'Upscaling' by using the surrogate model from a completed BOP-Elites run optimising a 25x25 solution archive and  predicting solutions for a 50x50 solution archive. All results in this section relate to Black-Box descriptor problems. Methods presented in the top of the table represent methods using the original budget (1250 evaluations), where both BOP-Elites and SPHEN use their surrogate methods to 'upscale', while those in the bottom rows are with an enlarged budget of 50k evaluations. 

Recall that 'upscaling' is neccesary because GP surrogate based methods cannot naively handle enough observations to fill a solution archive as large as 50x50. However, the results in Table \ref{tab:upscale} indicate that upscaling is a very effective method for extending surrogate based methods to high resolution problems. In each of the cases, both SPHEN and BOP-Elites provide competitive PMs, with BOP-Elites providing the best score of the PM methods. Interestingly, the Mishra Bird Function domain is again dominated by Sobol sampling methods for reasons discussed in the previous section.

\begin{figure*}[!t]
    \centering
    \subfloat[BOP-Elites, limited Budget]{\includegraphics[width=0.3\linewidth]{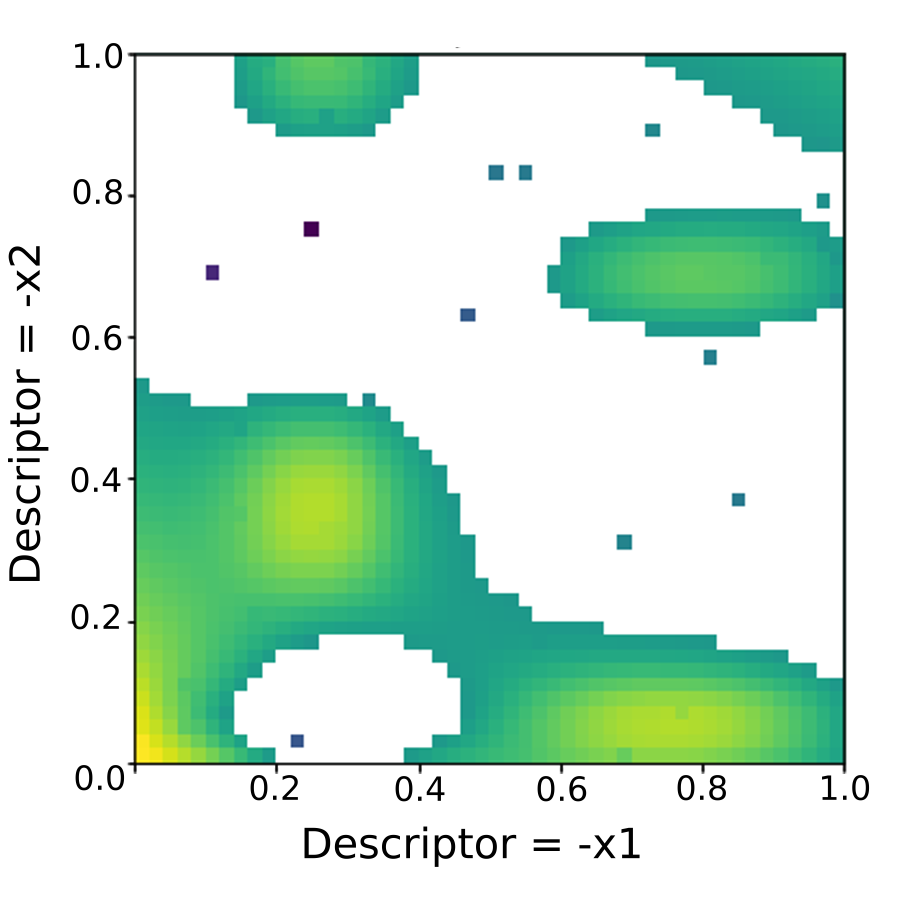}} 
    \subfloat[Sobol Sampling, 50k]{\includegraphics[width=0.3\linewidth]{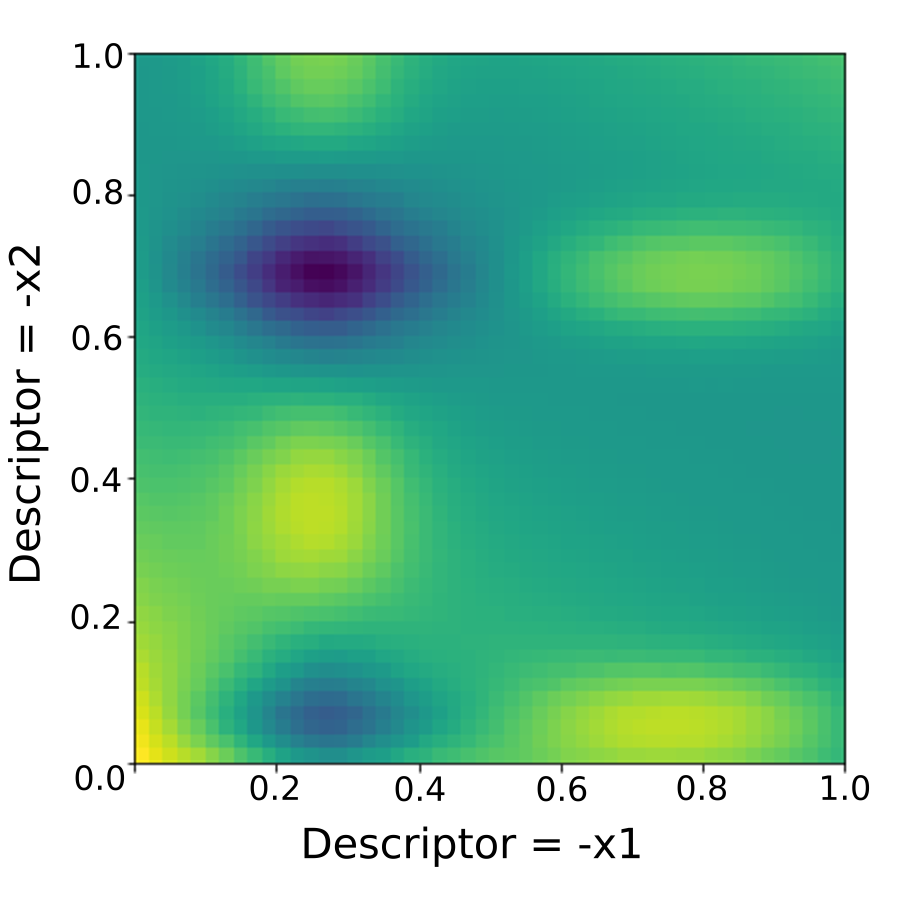}} 
    \subfloat[Upscaled BOP-Elites]{\includegraphics[width=0.3\linewidth]{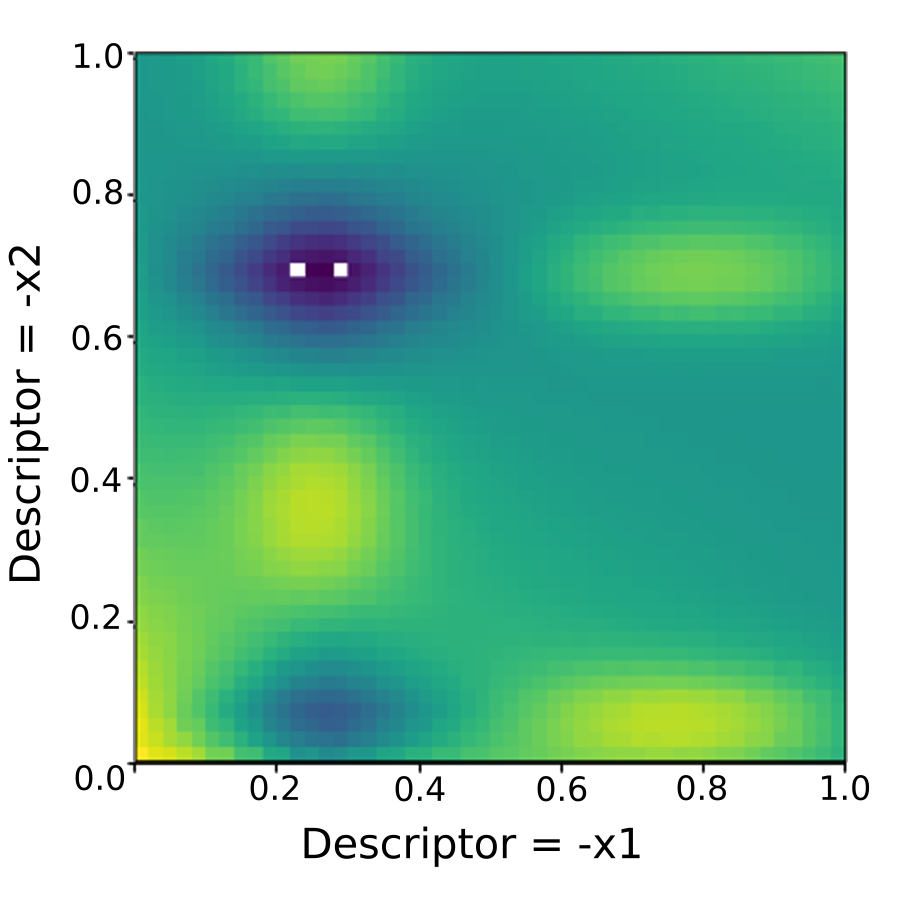}}
    \caption{Performance of a typical run of limited-budget BOP-Elites vs. Sobol sampling and upscaling on the Mishra's bird function.}
    \label{fig:Mishra-archives}
\end{figure*}

Figure \ref{fig:Mishra-archives} compares a $[50\times50]$  solution archive from an indicative run of BOP-Elites, a Sobol Sampling run with a 50k budget and an upscaled BOP-Elites solution. Fig. \ref{fig:Mishra-archives}(a) shows BOP-Elites is incapable of sampling each point of a 50x50 solution archive with its limited budget and therefore does not fill all the regions, but focuses on filling high performing areas first. In Fig. \ref{fig:Mishra-archives}(b) Sobol sampling with 50k evaluations performs well in filling the regions with high performing solutions, and in Fig. \ref{fig:Mishra-archives}(c) BOP-Elites Prediction Map upscaled from the 25x25 solution almost entirely fills the solution map trained on only 1.25k observations achieving 99\% of the best observed solution archive and outperforming MAP-Elites with 50k observations. 2 missing regions appear in areas of low performance

    

\begin{table}[!]
\begin{center}
\caption{A comparison of performance on multiple benchmarks given a $50\times50$ solution archive with the addition of an 'upscaled' $25\times25$ BOP-Elites PM }\label{tab:upscale}
\small
\setlength{\tabcolsep}{3pt}
\begin{tabular}{r@{\hspace{3pt}}rr@{\hspace{3pt}}rr@{\hspace{3pt}}rr}
  \hline
  QD score: 50x50& \multicolumn{2}{c}{Mishra} & \multicolumn{2}{c}{Robot Arm} & \multicolumn{2}{c}{Rosenbrock 6d}\\
  \cline{2-7}
  & Mean & S.E. & Mean & S.E.  & Mean & S.E.   \\
  \hline
  BOP-Elites         & 1.49e5 & 48.27 & 986.79 & 1.51  &  1.60e6 & 157.31 \\
  SAIL               & 1.30e5 & 771.27 & 897.66 & 10.07  &  1.18e6 & 1.31e4 \\ 
  SAIL - PM          & 3.07e5 & 22.71 & 1643.03 & 1.13 & 3.27e6 & 7.24e4 \\
  BOP-Elites - PM    & 3.05e5 & 38.99 & 1016.16 & 9.10 & 2.93e6 & 9.20e4 \\
  MAP-Elites         & 3.08e4 & 182.98 & 195.78 & 1.07  & 3.18e5 & 1.21e4   \\
  BOP-Elites-upscale  & \cellcolor{gray!25}3.08e5 & 51.99 & \cellcolor{gray!25}1764.56 & 5.51& \cellcolor{gray!25}3.33e6 & 3.10e4 \\
  \hline
  MAP-Elites 50k          & 3.02e5 & 152.1 & 1742.00 & 7.78 & 3.24e6 & 6.79e4 \\
  Sobol 50k            & \cellcolor{gray!25}3.10e5 & 8.05 & 267.07 & 1.15 & 1.75e6 & 618.10 \\
\hline
\end{tabular}
\end{center}
\end{table}


\section{Conclusion}

In this work, we proposed BOP-Elites, a Bayesian Optimisation algorithm for quality diversity search, and formulated a novel acquisition function, EJIE$^+$, which drastically outperforms MAP-Elites in both the black-box and white-box descriptor case using a form of Expected Improvement acquisition function for QDO problems and consistently outperforms both SAIL and SPHEN in the white-box and black-box optimisation tasks we explored. We note that BOP-Elites is extremely sample efficient, requiring far fewer observations than QD methods that do not use surrogate models and outperforming all surrogate assisted methods at the same budget when using the QD-score to assess the solution archive. BOP-Elites consistently produces solution archives with objective performance approaching optimality where the search budget is sufficient to fill the solution archives - our experiments suggest a budget of twice the number of regions is sufficient in many cases - and additionally produces surrogate models that can be leveraged to make high quality Prediction Maps for solution archives with more regions, a process we call upscaling. When the descriptor functions do not require modelling, as in the use-case for surrogate assisted algorithms (SAIL etc.), EJIE shows excellent performance, converging to near-optimality even faster. 

BOP-Elites is the first Bayesian Optimisation, rather than surrogate assisted, algorithm for QD problems and is able to model both the objective and descriptor values. Looking to the future, the BOP-Elites framework can be employed to test other acquisition functions such as Knowledge gradient, Entropy search etc. BOP-Elites can be generalised to working on noisy problems by implementation of an adapted Noisy-Expected-Improvement or Knowledge Gradient or other suitable acquisition functions. 

BOP-Elites outperforms our comparison algorithms in all of the cases presented. These results include the parsec domain which exhibits invalid regions and the robot arm problem which features unreachable regions. BOP-Elites provides a rapid improvement to its solution archive in the early stages, filling the archive with high quality points drawn from high-performing regions of the search space. This behaviour encourages good coverage, filling empty regions, starting with the high performing regions first. This is likely to be a desirable trait of the algorithm, as in
Fig.~\ref{fig:robotarm_convergence}a, early stopping of the algorithm, or a smaller budget will still produce exploration of the high performing regions.

Our results show that BOP-Elites produces surrogate models that can 'upscale', i.e., surrogates built for a solution archive with few regions can be used to generate a solution archive with many regions. This may be an interesting area of further development if we were able to populate a coarse solution archive with few observations and then actively refine the solution archive in areas of interest, e.g., subdivide regions with high density or high performance, spending more observations  in areas that prove the most useful. 

BOP-Elites slows down computationally in the prediction step when the Gaussian Process is built with many points, but this could be mitigated by the use of sparse GP techniques or alternative surrogate models. Gaussian Processes are generally considered to be ill-suited to high dimensional problems (beyond 10 dimensions) but BOP-Elites could also utilise other surrogate models to extend to much higher dimensions.

In addition to its strong performance, we consider BOP-Elites a proof of concept for continuing exploration of Bayesian Optimisation for QD problems and hope this work shows that, while QD research tends to come from the field of evolutionary computation, it is an interesting and rich area of future research in future applications for surrogate based methods and Bayesian Optimisation in particular.

\section*{ACKNOWLEDGEMENTS}
The first author would like to acknowledge funding from EPSRC through grant EP/L015374/1.

\bibliographystyle{IEEEtran}
\bibliography{ref,surrogate_me}

\end{document}